\magnification=\magstep1

\input amssym.def
\input amssym.tex
\font\twelvebf=cmbx12
%
%
\font\teneusm=eusm10
\font\seveneusm=eusm7
\font\fiveeusm=eusm5
\newfam\eusmfam
\textfont\eusmfam=\teneusm
\scriptfont\eusmfam=\seveneusm
\scriptscriptfont\eusmfam=\fiveeusm

\font\tenmib=cmmib10
\font\sevenmib=cmmib7
\font\fivemib=cmmib5
\newfam\mibfam
\textfont\mibfam=\tenmib
\scriptfont\mibfam=\sevenmib
\scriptscriptfont\mibfam=\fivemib

\font\tenss=cmss10
\font\sevenss=cmss8 scaled 833
\font\fivess=cmr5
\newfam\ssfam
\textfont\ssfam=\tenss
\scriptfont\ssfam=\sevenss
\scriptscriptfont\ssfam=\fivess
\def\ss{\fam\ssfam}
\thinmuskip = 2mu
\medmuskip = 2.5mu plus 1.5mu minus 2.1mu  
\thickmuskip = 4mu plus 6mu
\def\loosegraf#1\par{{%
\baselineskip=13.4pt plus 1pt \lineskiplimit=1pt \lineskip=1.3 pt
#1\par}}

\def\cc{{\Bbb C}}
\def\pp{{\Bbb P}}
\def\eps{\varepsilon}
\def\zz{{\Bbb Z}}
\def\rr{{\Bbb R}}
\def\ss{{\Bbb S}}
\def\bb{{\Bbb B}}

\def\calc{{\cal C}}


\centerline{\bf EXTENSION PROPERTIES OF MEROMORPHIC MAPPINGS WITH}

\medskip
\centerline{\bf VALUES IN NON-K\"AHLER COMPLEX MANIFOLDS}

\medskip
\centerline{\rm S. Ivashkovich\footnote{*}{This research was partially done 
during the author's stay at MSRI (supported in part by NSF grant DMS-9022140)
and at MPIM. I would like to give my thanks to both institutions for their 
hospitality. AMS subject 
classification: 32 D 15. Key words: meromorphic map, plurinegative metric,
Hartogs-type extension theorem, spherical shell, analytic cycle.}}


\bigskip
\noindent\bf
0. Introduction.

\smallskip\noindent\sl
0.1. Statement of the main result.
\smallskip\rm 
Denote by $\Delta (r)$ the disk of radius $r$ in 
$\cc $, $\Delta :=\Delta (1)$,  and for $0<r<1$ denote by  $A(r,1):= \Delta
\setminus \bar\Delta (r)$  an annulus in $\cc $. $\Delta^n(r)$ denotes 
the polydisk of radius $r$ in $\cc^n$ and  $\Delta^n:=\Delta^n(1)$.
Let $X$ be a compact complex manifold  and consider a
meromorphic mapping $f$ from the  ring domain $\Delta^n \times
A(r,1)$ into $X$. In this paper we shall study the following  

\smallskip\noindent
{\sl Question. \it Suppose
 we know that for some non-empty open subset $U\subset \Delta^n$ our map
 $f$ extends onto $U\times \Delta $. What is the maximal  $\hat U\supset U$ 
 such that $f$ extends meromorphically onto $\hat U\times \Delta $?}

\smallskip
 This is the so called Hartogs-type extension problem.  If $\hat U =\Delta^n$ 
 for any $f$ with values in our $X$ and any initial (non-empty!) $U$ then one says that the Hartogs-type 
 extension theorem holds for 
 meromorphic mappings into this $X$. 
 For $X=\cc $, i.e., for holomorphic functions the Hartogs-type extension theorem 
 was proved by F. Hartogs in [Ha]. If $X=\cc\pp^1$, i.e., for meromorphic 
 functions the result is due to E. Levi, see [Lv]. Since then the Hartogs-type 
 extension 
 theorem has been proved in at least two essentially more general cases than just 
 holomorphic or meromorphic functions. Namely for mappings into K\"ahler 
 manifolds and into  manifolds carrying  complete Hermitian metrics of 
 non-positive holomorphic sectional curvature, see [Gr], [Iv-3], [Si-2], [Sh-1].
 
 The goal of this paper is to initiate the systematic study of extension 
 properties of meromorphic mappings with values in non-K\"ahler complex 
 manifolds. Let $h$ be some Hermitian metric on a complex manifold $X$ and 
 let $\omega_h$ be the associated $(1,1)$-form. We call $\omega_h$ (and $h$ itself)
 {\it pluriclosed} or $dd^c$-closed  if $dd^c\omega_h=0$. In the sequel we shall not distinguish 
 between Hermitian metrics and their associated forms. The latters we shall call 
 simply metric forms. 
 
Let $A$ be a subset of $\Delta^{n+1}$ of Hausdorff $(2n-1)$-dimensional measure zero. 
Take 
a point $a\in A$ and a complex two-dimensional plane $P\ni a$ such that $P\cap A$ is  
of zero length. A sphere $\ss^3= \{ x\in P : \Vert x-a\Vert =\eps \} $ with $\eps $ 
small will be called a "transversal sphere" if in addition $\ss^3\cap A=\emptyset $. 
Take a non-empty open $U\subset \Delta^n$ and set $H^{n+1}_U(r)=
\Delta^n\times A(r,1)\cup U\times \Delta $. We call
 this set the Hartogs figure over $U$.

\break 
  \smallskip\noindent\bf
 Main Theorem. {\it Let $f:H^{n+1}_U(r)\to X$ be a meromorphic map into a compact
 complex manifold $X$, which admits a Hermitian metric $h$, such that the 
 associated $(1,1)$-form $\omega_h$ is $dd^c$-closed. Then $f$ extends to a
 meromorphic map  $\hat f:\Delta^{n+1}\setminus A   \to X$, where $A$
 is a complete $(n-1)$-polar, closed subset of $\Delta^{n+1}$ 
 of Hausdorff $(2n-1)$-dimensional measure zero. Moreover, if $A$ is the minimal 
 closed subset such that $f$ extends onto $\Delta^{n+1}\setminus A$ and 
 $A\not= \emptyset $,  then for every transversal 
 sphere $\ss^3\subset \Delta^{n+1}\setminus A$  its image $f(\ss^{3})$ is not  
 homologous to zero in $X$.

}

\smallskip\noindent\bf
Remarks 1. \rm A (two-dimensional) \it spherical shell \rm in a complex
manifold $X$ is the image $\Sigma $ of the standard sphere $\ss^3\subset
\cc^2$ under a holomorphic map of some neighborhood of $\ss^3$ into $X$
such that $\Sigma $ is not homologous to zero in $X$. Main Theorem states that if
the singularity set $A$ of our map $f$ is non-empty, then $X$ contains
spherical shells.

\smallskip\noindent\bf
2. \rm If again $A\not= \emptyset $ then, because 
$A\cap H^{n+1}_U(r)=\emptyset $ the restriction $\pi\mid_A:A\to \Delta^n$ of the 
natural projection 
$\pi :\Delta^{n+1}\to \Delta^n$ onto $A$ is proper. Therefore $\pi (A)$ is 
$(n-1)$-polar subset in $\Delta^n$ of zero $(2n-1)$-dimensional measure. So, 
returning to our question, we see 
that $\hat U$ is equal to $\Delta^n$ minus a "thin" set.

\smallskip We shall give a considerable number fo examples illustrating results of 
this paper. Let us mention few of them.
 
\smallskip\noindent\bf
Examples. 1. \rm Let $X$ be the  Hopf surface
 $X=(\cc^2\setminus \{ 0\} )/(z\sim 2z)$ and  $f:\cc^2\setminus \{ 0\}\to X$ 
 be the canonical projection. The $(1,1)$-form $\omega ={i\over 2}{dz_1\wedge 
 d\bar z_1+dz_2\wedge d\bar z_2\over \Vert z\Vert^2}$ is well defined  on
 $X$ and $dd^c\omega =0$. In this example one easily sees that $f$ is not extendable 
to zero and 
 that the image of the unit sphere from $\cc^2$ is not homologous to zero in 
 $X$. Note also that $dd^cf^*\omega = dd^c\omega = -c_4\delta_{\{ 0\} }dz\wedge d\bar z$,
where $c_4$ is the volume of the unit ball in $\cc^2$ and $\delta_{\{ 0\} }$ is 
the delta-fuction. 
 
\smallskip\noindent\bf
2. \rm In \S 3.6 we construct Example 3.7 of a $4$-dimensional compact complex 
manifold $X$ and a holomoprhic mapping $f:\bb^2\setminus \{ a_k\} \to X$, where 
$\{ a_k\} $ is a sequence of points converging to zero, such that $f$ cannot be 
meromorphically extended to the neighborhood of any $a_k$.

\smallskip\noindent\bf
3. \rm We also construct an Example 3.6 where the singularity set $A$ is of 
Cantor-type and pluripolar.

\smallskip This shows that the 
type of singularities described in our Main Theorem may occur. At the same time 
it should be noticed that we don't know if this $X$ can be endowed with a pluriclosed 
metric.

 \smallskip\noindent\bf
 4. \rm Consider now the Hopf three-fold $X=(\cc^3\setminus \{ 0\} )/(z\sim 2z)$.
 The analogous metric form $\omega ={i\over 2}{dz_1\wedge d\bar z_1+
 dz_2\wedge d\bar z_2+dz_3\wedge d\bar z_3\over \Vert z\Vert^2}$ is not longer 
 pluriclosed but only plurinegative (i.e. $dd^c\omega \le 0$). Moreover, if we 
 consider $\omega $ as a 
 bidimension $(2,2)$ current, then it will provide us a natural obstruction 
 for the existence of a  pluriclosed metric form on $X$. Natural projection 
 $f:\cc^3\setminus \{ 0\} \to X$ has singularity of codimension three and 
 $X$ doesn't contains spherical shells of dimension two (but contain a spherical
 shell of dimension three).

\smallskip\rm We also prove the Hartogs-type extension result for  mappings 
into (reduced, normal) complex spaces with $dd^c$-negative metric forms, see 
Theorem 2.2.
More examples, which are usefull for the understanding of the extension 
properties of meromorphic mappings into non-K\"ahler manifolds are given in 
the last paragraph. There also a general Conjecture is formulated.

\smallskip
 All compact complex surfaces admit pluriclosed Hermitian
metric forms. Therefore we have

\smallskip\noindent\bf
Corollary 1. {\it If $X$ is a compact complex surface, then: 

(a) every meromorphic
map $f:H^{n+1}_U(r)\to X$ extends onto $\Delta^{n+1}\setminus A$, where $A$ 
is an analytic set of pure codimension two;

(b) if $\Omega $ is a Stein surface and $K\Subset \Omega $ is a compact with 
connected complement, then every meromorphic map $f:\Omega \setminus K\to X$ 
extends onto $\Omega \setminus \{ ${\it finite set }$\} $. If this set is 
not empty (respectively, if $A$ from (a) is non-empty), then $X$ is of
class VII in the Enriques-Kodaira classification;

(c) if $f:\Omega \setminus K\to X$ is as in (b) but $\Omega $ of dimension 
at least three, then $f$ extends onto the whole $\Omega $.
}

\smallskip\noindent\bf
Remarks 1. \rm The fact that in the case of surfaces $A$ is a genuie analytic subset of pure 
codimension two requires some additional (not complicated) considerations and is given 
in \S 3.4, where also some other cases when $A$ can be proved to be analytic are 
discussed.

\smallskip\noindent\bf 
2. \rm
 A wide class of complex manifolds without
spherical shells is  for example the  class of such manifolds $X$ that the
Hurewicz homomorphism  $\pi_3(X)\to H_3(X,Z)$ is trivial.

\smallskip\noindent\bf
3. \rm Main Theorem  was proved in [Iv-2] under an additional
(very restrictive)
assumption: \it the manifold $X$ doesn't contain rational curves. \rm In
this case meromorphic maps into $X$ are  \it holomorphic \rm . Also in [Iv-2] nothing 
was proved about the structure of the singular set $A$. 

\smallskip\noindent\bf
4. \rm There is a hope that the surfaces with spherical shells could be 
 classified, as well as surfaces containing at least 
one rational curve. Therefore the following somewhat suprising 
speculation, which immediately follows from Corollary 1, could be of some 
interest:

\smallskip\noindent\bf
Corollary 2. {\it 
If a compact complex surface $X$ is not "among the known ones" then for every 
domain $D$ in a Stein surface every meromorphic mapping $f:D \to X$ 
is in fact holomorphic and extends as a holomorphic mapping $\hat f:\hat D\to X$
of the envelope of holomorphy $\hat D$ of $D$ into $X$.
}

\smallskip\rm At this point let us note that the notion of a sphericall shell, 
as we understand it here, is different from the notion of {\it global 
spherical shell} from [Ka-1].
 
\smallskip\noindent\bf
5. \rm A real two-form $\omega $ on a complex manifold $X$ is said to
"tame" the complex structure $J$ if for any non-zero tangent vector $v\in
TX$ we have $\omega (v,Jv)>0$. This is equivalent to the property that the 
$(1,1)$-component
$\omega^{1,1}$ of $\omega$ is strictly positive. Complex manifolds admitting a
\it closed \rm form, which tames the complex structure, are of special interest.
The class of such manifolds  contains all K\"ahler manifolds. On the
other hand, such metric forms are $dd^c$-closed. Indeed, if
$\omega =\omega^{2,0} +\omega^{1,1} +\bar \omega^{2,0}$ and $d\omega=0$, then
$\partial \omega^{1,1}=-\bar \partial \omega^{2,0}$.
Therefore $dd^c\omega^{1,1}=2i\partial \bar \partial \omega^{1,1}=0$. So the Main Theorem 
applies to meromorphic mappings into such manifolds. In fact, the technique of the
proof gives more:

\smallskip\noindent\bf
Corollary 3. \it Suppose that a compact complex manifold  $X$ admits a
strictly positive $(1,1)$-form, which is the $(1,1)$-component of a closed
form. Then every meromorphic map 

\noindent $f:H^{n+1}_U(r)\to X$ extends onto $\Delta^{n+1}$.

\smallskip\rm
This statement generalizes the Hartogs-type extension theorem for 
meromorphic mappings into K\"ahler manifolds from [Iv-3], but this 
generalization cannot be obtained by the methods of [Iv-3] and result from [Si-2]
involved there. The reason is simply that the upper levels of 
Lelong numbers of pluriclosed (i.e., $dd^c$-closed) currents are no longer 
analytic (also integration 
by parts for $dd^c$-closed forms doesn't work as well as for $d$-closed ones).

\smallskip
It is also natural to consider  the extension of meromorphic mappings from
\it singular \rm spaces. This is equivalent to considering  multivalued
meromorphic correspondences from smooth domains, and this reduces to singlevalued
maps into symmetric powers of the image space, see \S 3 for details. However,
one pays a price for these reductions. In this direction we construct in \S 3
{\sl Example 3.5}, which shows that a manifold possessing the Hartogs extension 
property for singlevalued mappings may not  possess it for 
multivalued ones. The reason is  that ${\sl Sym}^2(X)$ may contain a spherical shell, if 
even $X$ contains none.

\smallskip\noindent\sl
0.2. Sketch of the proof.
\smallskip\rm
Let us give a brief outline of the proof of the main theorem. We first  
consider the case of dimension two, i.e., $n=1$. For $z\in \Delta $ set $\Delta_z:=
\{ z\} \times \Delta $. For a meromorphic map $f:H^2_U(r)\to (X,\omega )$ denote by 
$a(z)=area_{\omega }f(\Delta_z)=\int_{\Delta }f|_{\Delta_z}^*\omega  $ - the area of the 
image of the disk $\Delta_z$. This is well defined for $z\in U$ after shrinking 
$A(r,1)$ if necessary.

\smallskip\noindent\sl
Step 1. \rm Using $dd^c$-closedness of $\omega $ (and therefore of $f^*\omega $) we 
show that for "almost every" sequence $\{ z_n\} \subset U$ converging to the boundary, 
areas $a(z_n)$ are uniformly bounded and converge to the area of $f(\Delta_{z_{\infty }})$, 
here $z_{\infty }\in \partial U\cap \Delta $ is the limit of $\{ z_n\} $. This means in 
particular that $f_{z_{\infty }}:=f|_{\{ z_{\infty }\} \times A(r,1)}$ extends onto 
$\Delta_{z_{\infty }}$. And then we show that $f$ can be extended holomorphically onto 
$V\times \Delta $, where $V$ is a neighborhood of $z_{\infty }$. Therefore if $\hat U$ is 
the maximal open set s.t. $f$ can be extended onto $H^2_{\hat U}(r)$, then $\partial \hat U
\cap \Delta $ should be "small". In fact we show that  $\partial \hat U\cap \Delta $  is 
of harmonic measure zero, see Lemmas 2.3, 2.4.

\smallskip\noindent\sl
Step 2. \rm Interchanging coordinates in $\cc^2$ and repeating the Step 1 we see that $f$ 
holomorphically extends onto $\Delta^2\setminus (S_1\times S_2)$, where $S_1$ and $S_2$ are 
compacts (after shrinking) of harmonic measure zero. We can use shrinkung here, because 
subsets of harmonic measure zero in $\cc $ are of Hausdorff dimension zero. Set $S=S_1\times 
S_2$. Smooth form $T:=f^*\omega $ on $\Delta^2\setminus S$ has coefficients in $L^1_{loc}(\Delta^2)$ 
and therefore has trivial extension $\tilde T$ onto $\Delta^2$, see Lemma 3.3 form [Iv-2] and 
Lemma 2.1. We prove that $\mu :=dd^c\tilde T$ is a non-positive measure supported on $S$.

\smallskip\noindent\sl
Step 3. \rm Take a point $s_0\in S$ and, using the fact that $S$ is of Hausdorff dimension zero, 
take a small ball $B$ centered at $s_0$ such that $\partial B\cap S=\emptyset $. Now we have
two possibilities. First: $f(\partial B)$ is not homologous to zero in $X$. Then 
$\partial B$ represent a spherical shell in $X$, as it was said in the remark after 
the Main Theorem. Second: $f(\partial B)\sim 0$ in $X$. Then we can prove, see Lemmas 2.5,
2.8 that $\tilde T$ is $dd^c$-closed and 
consequently it be written in the form $\tilde T=i(\partial\bar\gamma - 
\bar\partial\gamma )$, where $\gamma $ is some $(0,1)$-current on $B$, which is 
smooth on $B\setminus S$. This allows to estimate the area function $a(z)$ in the 
neighborhood of $s_0$ and extend $f$.

\medskip\noindent\sl
Step 4. We consider now the case $n\ge 2$. 
\smallskip\rm Using case $n=1$ by sections we extend $f$ onto $\Delta^{n+1}\setminus A$ 
where $A$ is complete pluripolar of Hausdorff codimension four. Then take a transversal 
to $A$  at point $a\in A$ complex two-dimensional direction and decompose the neighborhood $W$ 
of $a$ as $W=B^{n-1}\times B^2$, where $A\cap (B^{n-1}\times \partial B^2)=\emptyset $. If $f(\{ a\} \times 
\partial B^2)$ is homologous to zero then we can repeat Step 3 "with parameters". This will give 
a uniform bound of the volume of the two-dimensional sections of the graph of $f$. Now we 
are in position ty apply the Lemma 1.3 (which is another main ingredient of this paper) 
to extend $f$ onto $W$.

\smallskip\noindent\bf
Remark. \rm We want to finish this {\it Introduction} with a brief account of 
existing methods of extension of meromorphic mappings. The {\sl First} method, 
based on Bishop's extension theorem for analytic sets (appearing here as the graphs 
of mappings) and clever integration by parts was introduced by P. Griffiths 
in [Gr], developed by B. Shiffman in [Sh-2] and substantially 
enforced by Y.-T. Siu in [Si-2] (where the  Thullen-type extension theorem is proved 
for mappings into K\"ahler manifolds), using 
his celebrated result on the 
analyticity of upper level sets of Lelong numbers of closed positive currents.
The latter was by the way inspired by the extension theorem just mentioned. 
Finally, in [Iv-3] the Hartogs-type extendibility for the mappings into 
K\"ahler manifolds was proved using the result of Siu and a somewhat 
generalized classical method of "analytic disks". This method works well for 
mappings into K\"ahler manifolds.

The {\it Second} method, based on the Hironaka imbedded resolution of singularities 
and lower estimates of Lelong numbers was proposed in [Iv-4] together with 
an example showing the principal difference between K\"ahler and non-K\"ahler cases. 
This method implies  the Main Theorem of this paper for $n=1,2$ 
(this was not stated in [Iv-4]). 
However, further increasing of $n$ meets technical difficulties at 
least on the level of the full and detailed proof of Hironaka's theorem 
(plus it should be accomplished with the detailed lower estimates of the Lelong 
numbers by blowings-up).

The {\it Third } method is therefore proposed  in this paper and is based on the 
Barlet cycle space theory. It gives definitely stronger and more 
general  results than the previous two  and is basically {\it much more simple}.
Key point is Lemma 1.3 from \S 1.
An important ingredient of last two methods is the notion of meromorphic 
family of analytic subsets and especially {\sl Lemma 2.4.1} from [Iv-4] about such 
families. The reader is therefore supposed  to be familiar with \S \S 2.3 and 
2.4 of [Iv-4] while reading  proofs of both Lemma 1.3 and  Main Theorem.

\smallskip\rm 

I would like to give my thanks to the refery, who pointed to me a gap in the proof of 
the analyticity of the singular set.    

\bigskip
\def\longpoints{\leaders\hbox to 0.5em{\hss.\hss}\hfill \hskip0pt}

\smallskip\noindent\sl
\smallskip\rm
\vbox{\hsize=4.5truein
\baselineskip=10pt \sevenrm
\bigskip\noindent\bf
Table of content\medskip \sl
\noindent
\line{0. Introduction. \longpoints pp. 1 - 5}
\medskip\noindent
\line{1. Meromorphic mappings and cycle spaces. \longpoints pp. 5 - 14}
\medskip\noindent

\line{2. Hartogs-type extension and spherical shells. \longpoints pp. 14 - 25}
\medskip\noindent
\line{3. Examples and open questions. \longpoints pp. 25 -  34}
\medskip\noindent
\line{ References. \longpoints pp. 34 - 37}
}
\bigskip

\smallskip\noindent\bf
1. Meromorphic mappings and cycle spaces.
\smallskip\noindent\sl
1.1. Cycle space associated to a meromorphic map.
\smallskip\rm
We shall freely use the results from the theory of cycle spaces
developed by D.Barlet, see [Ba-1]. For the English spelling of
Barlet's terminology we refer to [Fj].
Recall that an analytic $k$-cycle 
in a complex space $Y$ is a formal sum $Z=\sum_jn_jZ_j$, where $\{ Z_j\} $ is a 
locally finite sequence of  analytic
subsets (always of pure dimension
$k$) and $n_j$ are positive integers called multiplicities of $Z_j$.
Let $\vert Z\vert :=\bigcup_jZ_j$ be the support of $Z$. All complex spaces in this
paper are reduced, normal and countable at infinity. All cycles, if the 
opposite is not stated, are supposed to have  {\it connected} support. Set 
$A^k(r,1)=\Delta^k\setminus \bar\Delta^k(r)$.

Let $X$ be a normal, reduced complex space equipped with some Hermitian metric. 
Let a holomorphic mapping $f:\bar\Delta^n\times \bar A^k(r,1)\to X$ be given. 
We shall start with the following space
of cycles related to $f$. Fix some positive constant $C$ and consider the set 
${\cal C}_{f,C}$ of all
analytic $k$-cycles $Z$ in $Y:=\Delta
^{n+k}\times X$ such that:

(a) $Z\cap [\Delta^n\times \bar A^k(r,1)\times  X] = \Gamma_{f_{z}}\cap
[\bar A^k_z(r,1)\times X]$ for
some $z\in \Delta^n$, where $\Gamma_{f_z}$ is the graph of the restriction 
$f_z:=f\mid_{A^k_z(r,1)}$. Here $A^k_z(r,1):=\{ z\} \times A^k(r,1)$. This 
means, in particular, that for this $z$ the
mapping $f_{z}$ extends meromorphically from $\bar A_{z}^k(r,1)$ onto
$\bar\Delta_{z}^k:=\{ z\} \times \bar\Delta^k$.

(b) ${\sl vol}(Z)< C$ and the support $\vert Z\vert $ of $Z$ 
is connected.

\smallskip
We put $\calc_{f}:=\bigcup_{C>0}\calc_{f,C}$ and we are going to show that 
${\cal C}_{f}$ is an analytic space of finite dimension in a neighborhood of 
each of its points.

\smallskip
Let $Z$ be an analytic cycle of dimension $k$ in a (reduced, normal) complex
space $Y$. In our applications $Y$ will be $\Delta^{n+k}\times X$. By a
coordinate chart  adapted to $Z$ we shall understand an  open set $V$
 in $Y$ such that $V\cap \vert Z\vert \not=\emptyset $ together with an
isomorphism $j $ of $V$ onto a closed subvariety $\tilde V$ in the
neighborhood
 of $\bar\Delta^k\times \bar\Delta^q$ such that $j^{-1}(\bar\Delta^k
 \times
\partial \Delta^q)\cap \vert Z\vert =\emptyset $. We shall denote such a chart
by $(V,j )$. The image $j (Z)$ of cycle $Z$ under isomorphism $j $
is the image of the underlying analytic set
together with multiplicities. Sometimes we shall, following Barlet, denote: 
$\Delta^k=U, \Delta^q=B$ and call the quadruple $E=(V,j,U,B)$ a  \it 
scale \rm adapted to $Z$.

If ${\sl pr} :\cc^k \times \cc^q\to \cc^k $ is the natural
projection, then the restriction ${\sl pr}\mid_{j(Z)}:j(Z)\to \Delta^k $ is a 
branched 
covering of degree
say $d$. The number $q$ depends on the imbedding dimension of $Y$ (or $X$ in our
case). Sometimes we shall skip  $j$ in our notations. The branched
covering
${\sl pr}\mid_Z:Z\cap (\Delta^k\times \Delta^q)\to
\Delta^k $ defines in a natural way a mapping $\phi_Z :\Delta^k\to {\sl Sym}^d
(\Delta^q)$ - the $d$-th 
 symmetric power of $\Delta^q$ - by setting  
 $\phi_Z (z)=({\sl pr}\mid_Z)^{-1}(z)$. This
allows to represent a cycle $Z\cap \Delta^{k+q}$ with $\vert Z\vert \cap
(\bar\Delta^k\times \partial\Delta^q)=\emptyset $ as the  graph of a $d$-valued
holomorphic map.

\smallskip\rm
Without loss of generality we suppose that our holomorphic mapping  $f$ is
defined
on $\Delta^n(a)\times A^k(r_1,b)$ with $a,b>1,r_1<r$. Now, each
$Z\in {\cal C}_{f}$ can be covered by a \it finite \rm number of adapted
neighborhoods $(V_{\alpha },j_{\alpha })$. Such covering we be called an 
adapted covering. Denote the union $\bigcup_{\alpha }V_{\alpha }$
by $W_Z$. Taking this covering $\{ (V_{\alpha },j_{\alpha })\{ $
to be small enough, we can  further suppose that:

\smallskip\noindent
{\it (c)  if $V_{\alpha_1}\cap V_{\alpha_2}\not=\emptyset $, then on every
irreducible component of the intersection $Z\cap V_{\alpha_1}\cap V_{\alpha
_2}$ a point $x_1$ is fixed so that:

$(c_1)$ either there exists a polycylindric 
neighborhood $\Delta_1^k\subset \Delta^k$ of $pr(j_{\alpha_1}(x_1))$ such that
the
chart $V_{12}=j_{\alpha_1}^{-1}(\Delta_1^k\times \Delta^q)$ is
adapted to $Z$ and is
contained in $V_{\alpha_2}$, where $V_{12}$ is given the same imbedding
$j_{\alpha_1}$,

$(c_2)$ or this is fulfilled for $V_{\alpha_2}$ instead of $V_{\alpha_1}$;

\noindent
(d) if $V_{\alpha }\ni y$ with $p(y)\in \bar\Delta^n(c)\times A^k({r+1
\over 2},1)$, then $p(\bar V_{\alpha })\subset \bar\Delta^n({c+1\over 2})
\times A^k(r,1)$.
}

\smallskip
Here we denote by $p :\Delta^{n+k}\times X\to \Delta^{n+k}$ the 
 natural projection. Case $(c_1)$ can be  realized when the imbedding dimension of
$V_{\alpha_1}$ is smaller or equal to that of $V_{\alpha_2}$, and $(c_2)$
in the opposite case, see [Ba-1], pp. 91-92.

Let $E=(V,j,U,B)$ be a scale on the complex space $Y$. Denote by 
$H_Y(\bar U,{\sl sym}^d(B)):={\sl Hol}_Y(\bar U,{\sl sym}^d(B))$  the 
 Banach analytic set
of all $d$-sheeted analytic subsets on $\bar U\times B$, contained
in $j(Y)$. The subsets $W_Z$ together with the 
topology of uniform convergence on $H_Y(\bar U,{\sl sym}^d(B))$
define a (metrizable) topology on our cycle space ${\cal C}_{f}$, which
is equivalent to the topology of currents, see [Fj],[H-S].

 We refer the  reader to [Ba-1] for the definition of the isotropicity of the 
family of elements from $H_Y(\bar U,{\sl sym}^d(B))$ parametrized by some 
Banach analytic set ${\cal S}$.  Space $H_Y(\bar U,{\sl sym}^d(B))$ can be endowed by another 
(more rich) analytic structure. This new analytic space will be denoted by 
$\hat H_Y(\bar U,{\sl sym}^d(B))$.  The crucial property of this new structure 
is that the tautological family 
$\hat H_Y(\bar U,{\sl sym}^d(B))\times U'\to {\sl sym}^d(B)$ is isotropic 
in $H_Y(\bar U',{\sl sym}^d(B))$ for any 
relatively compact polydisk $U'\Subset U$, see [Ba-1].
In fact for isotropic families $\{ Z_s:s\in {\cal S} \} $ parametrized by 
Banach analytic sets the following  Projection Changing theorem of 
Barlet holds.

\smallskip\noindent\bf
Theorem (Barlet). {\it If the family $\{ Z_s:s\in {\cal S} \} \subset 
H_Y(\bar U,{\sl sym}^d(B))$ is isotropic, then for any scale 
$E_1=(V_1,j_1,U_1,B_1)$ in $U\times B$ adapted to some $Z_{s_0}$, there exists a 
neighborhood $U_{s_0}$ of $s_0$ in ${\cal S}$ such that $\{ Z_s:s\in U_{s_0}\} $ is 
again isotropic in $V_1$.
}
\smallskip\rm This means, in particular, that the mapping $s\to Z_s\cap V_1\subset 
H_Y(\bar U_1,{\sl sym}^d(B_1))$ is analytic, i.e., can be extended to a neighborhood
of any $s\in U_{s_0}$. Neighborhood means here a neighborhood in some complex 
Banach space where ${\cal S}$ is defined as an analytic subset.

This leads naturally to the following

\smallskip\noindent\bf
Definition 1.1. {\it A family ${\cal Z}$ of analytic cycles in an open set
$W\subset Y$, parametrized by a Banach analytic space ${\cal S}$, is called 
analytic
in a neighborhood of $s_0\in {\cal S}$ if for any scale $E$ 
adapted to $Z_{s_0}$
there exists a neighborhood $U\ni s_0$ such that  the family $\{ {\cal Z}_s:s\in U\}$ 
is isotropic.
}

\medskip
\noindent
{\sl 1.2. Analyticity of ${\cal C}_{f}$ and construction of ${\cal G}_f$ }.

\smallskip\rm
Let $f:\bar\Delta^n\times \bar A^k(r,1)\to X$ be our map.  Take a cycle  
$Z\in {\cal C}_{f}$ and a finite covering $(V_{\alpha },j_{\alpha })$ satisfying 
conditions (c) and (d). As above, put $W_Z=\bigcup V_{\alpha }$.
We want to show now that ${\cal C}_{f}$ is an analytic space of finite
dimension in a neighborhood of $Z$. We divide $V_{\alpha }$'s into two types.
\smallskip
\noindent\it 
Type 1. These are $V_{\alpha }$ as in (d). For them  put

$$
H_{\alpha }:=\bigcup_z
\{ [\Gamma_{f_z}\cap \bar A^k_z(r,1)\times X]\cap V_{\alpha }\}
\subset H_Y(\bar U_{\alpha }, {\sl Sym}^{d_{\alpha }}(B_{\alpha })).\eqno(1.2.1)
$$
\noindent\rm
 The union is
taken over all $z\in \Delta^n$ such that $V_{\alpha }$ is adapted to
$\Gamma_{f_z}$.

\noindent\it
Type 2. These are all others. For  $V_{\alpha }$ of this type we put $H_{\alpha }:=\hat H_Y
(\bar U_{\alpha },{\sl Sym}^{d_{\alpha }}(B_{\alpha }))$.

\smallskip\rm
All $H_{\alpha }$ are open sets in complex Banach analytic subsets and for $V_{\alpha }
$ of the first type they are of dimension $n$ and smooth. The latter follows from 
Barlet-Mazet theorem, which tells that if $h:A\to {\cal S}$ is a holomorphic 
injection of a finite dimensional analytic set $A$ into a Banach analytic 
set ${\cal S}$, then $h(A)$ is also an Banach analytic set of finite dimension, see [Mz].

For every irreducible component of $V_{\alpha }\cap V_{\beta }\cap Z_l$ we
fix a  point
$x_{\alpha \beta l}$ on this component (the subscript $l$ indicates the component),
and a 
chart $V_{\alpha }\cap V_{\beta }\supset (V_{\alpha \beta l},\phi_{\alpha
\beta l})\ni x_{\alpha \beta l}$ adapted to this
component as in (c). Put $H_{\alpha \beta l}:=\hat H(\Delta^k,
{\sl Sym}^{d_{\alpha \beta l}}(\Delta^p))$. In the sequel it will be convenient to introduce
an order on our finite covering $\{ V_{\alpha }\} $ and write $\{ V_{\alpha }\}_{\alpha =1}^N$.

Consider finite products $\Pi_{(\alpha )}H_{\alpha }$ and $\Pi_{(\alpha
\beta l)}H_{\alpha \beta l}$. In the second
product we take only triples with $\alpha <\beta $. These are Banach analytic
spaces and by
 the Projection Changing Theorem of Barlet, for each pair $\alpha <\beta $ we
have two
holomorphic mappings $\Phi_{\alpha \beta }:H_{\alpha }\to \Pi_{(l)}H_
{(\alpha \beta l)}$
 and $\Psi_{\alpha \beta }:H_{\beta }\to
\Pi_{(l)}H_{\alpha \beta l}$. This defines two holomorphic maps $\Phi ,\Psi
:\Pi_{(\alpha )}
H_{\alpha }\to \Pi_{\alpha <\beta ,l}H_{\alpha \beta l}$. The kernel 
${\cal A}$
 of this pair, i.e., the set of 
$h=\{ h_\alpha \} $
 with $\Phi (h)=\Psi (h)$, consists exactly from analytic cycles in the 
 neighborhood $W_{Z}$ of $Z$. This kernel is a Banach analytic set, and
 moreover the family ${\cal A}$ is an analytic family in $W_{Z}$ in the 
 sense of {\sl Definition 1.1}.
 \smallskip\noindent\bf
 Lemma 1.1. \it ${\cal A}$ is of finite dimension.
 \smallskip\noindent\bf
 Proof. \rm Take a smaller covering $\{ V_{\alpha }^{'},
 j_\alpha\} $ of
 $Z$. Namely, $V_{\alpha }^{'}=V_{\alpha } $ for  $V_{\alpha }$ of the
 first type and $V_{\alpha }^{'}=
 j_{\alpha }^{-1}(\Delta_{1-\eps }\times \Delta^p)$ for the second. In
 the same  manner define $H_{\alpha }^{'}$ and $H^{'}:=\Pi_{\alpha }H_{\alpha }^
 {'}$. Repeating the same construction as above we obtain a Banach analytic 
 set ${\cal A}^{'}$. We have a holomorphic
 mapping $K:{\cal A}\to {\cal A}^{'}$ defined by the restrictions.
 The differential $dK\equiv K$ of this map is a compact operator. 
 
 Let us show that we also have 
 an  inverse {\it analytic} map $F:{\cal A}^{'}\to {\cal A}$. The analyticity  
 of $F$ means more precisely that it should be defined in some neighborhood 
 of ${\cal A}^{'}$ in $H^{'}$. For scales $E_{\alpha } =(V_{\alpha },U_{\alpha }
 ,B_{\alpha },j_{\alpha })$ of the second type the mapping 
 $F_{\alpha }:{\cal A}
 ^{'}\to H_Y(\bar U_{\alpha },{\sl Sym}^{k_{\alpha }}B_{\alpha })$ is defined 
 by the isotropicity of the family ${\cal A}^{'}$ as in [Ba-1]. In particular, 
 this $F_{\alpha }$ extends analytically to a neighborhood in $H^{'}$(!) of 
 each point of ${\cal A}^{'}$.
 
 For scales $E_{\alpha } =(V_{\alpha },U_{\alpha }=U_{\alpha }^{'},
 B_{\alpha },j_{\alpha })$ of the first type define $F_{\alpha }$ as follows.
 Let $Y=(Y_{\alpha })$ be some point from $H^{'}$. Using the fact that 
 $H_{\alpha }=H_{\alpha }^{'}$ in this case, we can correctly define 
 $F_{\alpha }(Y):=Y_{\alpha }$ viewed as an element  of $H_{\alpha }$. This 
 directly defines $F_{\alpha }$ on the whole $H^{'}$. Analyticity is also 
 obvious. 
 
 Put $F:=\Pi_{\alpha }F_{\alpha }:{\cal A}^{'}\to {\cal A}$. $F$ is defined and 
 analytic in a neighborhood of each point of ${\cal A}^{'}$. 
 Observe further that  $id-dK\circ dF$ is Fredholm. Since  
 ${\cal A}^{'}\subset \{ h\in
 \Pi_{(i)}H_i^{'}: (id-K\circ F)(h)=0\} $, we obtain that ${\cal A}^{'}$ is
 an analytic subset in a complex manifold of finite dimension.
 \smallskip
 \hfill{q.e.d.}
 \smallskip
Therefore ${\cal C}_f$ is an analytic space of finite dimension in a neighborhood 
of each of its points. ${\cal C}_{f,C}$ are  open subsets of ${\cal C}_f$. 
Note further that for $C_1<C_2$ the set ${\cal C}_{f,C_1}$  is an open subset of 
${\cal C}_{f,C_2}$. This implies that for each irreducible component 
${\cal K}_C
$ of ${\cal C}_{f,C}$ there is a unique irreducible component ${\cal K}$ of 
${\cal C}_{f}$ containing ${\cal K}_C$ and moreover ${\cal K}_C$ is an 
open subset of ${\cal K}$. Of course, in general the dimension of irreducible 
components of ${\cal C}_f$ is not bounded, and in fact the space ${\cal C}_f$ 
is too big. Let us denote by ${\cal G}_f$ the union of  irreducible 
components of
${\cal C}_f$ that contain at least one irreducible cycle or, in other words, a 
cycle of the form $\Gamma_{f_z}$ for some $z\in \Delta^n$. 

Denote by 
${\cal Z}_f:=\{ Z_a:a\in {\cal C}_f\} $ the universal family.
In the sequel ${\cal B}_k(X)$ will denote  the Barlet space of 
compact analytic $k$-cycles in normal, reduced complex space $X$.
\smallskip\noindent\bf
Lemma 1.2. {\it 1. Irreducible cycles form an open dense subset ${\cal G}_f^0$ 
in ${\cal G}_f$.

2. The dimension of ${\cal G}_f$ is not greater than $n$.

3. If $k=1$, then all compact irreducible  components of cycles in ${\cal G}_f$ are rational.

}
\smallskip\noindent\bf
Proof. \rm 1. ${\cal G}_f^0$ is clearly open, this follows immediately from 
(4) and (6) of {\sl Lemma 2.3.1} in [Iv-4]. Denote by $\hat {\cal C}_f$ the 
normalization of ${\cal C}_f$ and denote by $\hat {\cal Z}_f$  the pull-back 
of the universal family under the normalization map ${\cal N}:\hat {\cal C}_f 
\to {\cal C}_f$. Consider the following "forgetting of extra compact components"
mapping $\Pi :\hat {\cal C}_f \to \hat {\cal C}_f$. Note that each cycle $Z\in 
\hat {\cal C}_f$ can be uniquely represented as $Z = \Gamma_{f_s} + \Sigma_{j=1}
^NB_s^j$, where each $B_s^j$  is a compact analytic $k$-cycle  
in $\Delta^k_s(r)\times X$ with connected support. Mark those $B_s^j$, which 
possess the following property: there is a neighborhood in $V$ of $Z$ in 
$\hat {\cal C}_f$ such that every cycle $Z_1\in V$ decomposes as $Z_1=\hat Z_1 +
B_1$, where $B_1$ is a compact cycle in a neighborhood of 
$B_s^j$ in the Barlet space ${\cal B}_k(X)$. Our mapping $\Pi :  \hat {\cal C}_f
\to \hat {\cal C}_f$ sends each cycle $Z$ to the cycle obtained from this $Z$ 
by deleting all the marked components. This is clearly an analytic map. Every 
irreducible cycle is clearly a fixed point of $\Pi $. Thus the set of fixed 
points is open in $\hat {\cal G}_f\subset \hat {\cal C}_f$ and so contains the whole 
$\hat {\cal G}_f$. 

Now we shall prove that every fixed point $Z$ of $\Pi $ is a limit of 
irreducible cycles. For the sequel remark that the compositions 
$\psi :=p\circ {\bf ev} : {\cal Z}_f \to \Delta^{n+k} $ and $\phi := p_1\circ 
{\bf ev} \circ \pi^{-1} : {\cal C}_f\to \Delta^n$ are well defined. Here 
$p_1:\Delta^{n+k}\times X\to \Delta^n$ is one more natural projection and 
${\bf ev} : {\cal Z}_f \to \Delta^{n+k}\times X $ is the natural evaluation map. 
Let $\phi (Z)=s\in \Delta^n$ and 
$Z = \Gamma_{f_s} + \Sigma_{j=1}^N B_s^j$. $Z$ being a fixed point of $\Pi $ 
means that in any neighborhood of $Z$ one can find a cycle $Z_1$ such that 
$Z_1=\Gamma_{f_{s_1}} + \Sigma_{j=2}^NB^j_{s_1}$, where $B_{s_1}^j$ are compact 
cycles close to $B^j_s$. Observe that every cycle in a neighborhood of 
$Z_1$ has the same form, i.e.,  in its decomposition $j\ge 2$, which follows 
from {\sl Lemma 2.3.1} from [Iv-4]. Since $Z_1$ is also a fixed point for 
$\Pi $, we can repeat this procedure $N$ times to obtain finally an irreducible 
cycle in a given neighborhood of $Z$.

We conclude that ${\cal G}_f^0$ is dense in  ${\cal G}_f$. 

\smallskip\noindent
2. Take an irreducible $Z\in {\cal G}_f^0\cap {\sl Reg}({\cal G}_f )$. Take a 
neighborhood $Z\in V\subset {\sl Reg}({\cal G}_f )$ that consists from 
irreducible cycles only. Then $\phi\mid_V:V\to \Delta^n$ is injective and 
holomorphic. Thus ${\sl dim}{\cal G}_f \le n$.

\smallskip\noindent
3. This part follows from {\sl Lemma 7} in [Iv-5] because every cycle from 
 ${\cal G}_f$ is a limit of analytic disks.
 
 \smallskip
 \hfill{q.e.d.}

 \smallskip\noindent\bf
 Definition 1.2. {\it  We shall call the space ${\cal G}_f$ the {\sl cycle space} 
 associated to a meromorphic map $f$.
 }
  
 \smallskip\rm
 Denote by ${\cal G}_{f,C}$  the open subset of ${\cal G}_f$ consisting 
 of $Z$ with ${\sl vol}(Z) < C$.

 \smallskip\noindent\sl
 1.3. Proof of the Main Lemma.
 \smallskip\rm
 Now we are ready to state and prove the main lemma of this paragraph, i.e. Lemma 1.3. 
 From now on we restrict 
 our universal family ${\cal Z}_{f}$ onto  ${\cal G}_{f}$ withought changing 
 notations. I.e. now ${\cal Z}_{f,C}:=\{ Z_a:a\in {\cal G}_{f,C}\} $, 
 ${\cal Z}_{f}:=\bigcup_{C>0}{\cal Z}_{f,C}$ and $\pi :{\cal Z}_{f}\to 
 {\cal G}_{f}$ is the natural projection. ${\cal Z}_{f}$ is a 
 complex space of finite dimension. We have an evaluation map
 $$
 {\bf ev}:{\cal Z}_f\to \Delta^{n+k}\times X,\eqno(1.3.1)
 $$
 \noindent
 defined by $Z_a\in {\cal Z}_f\to Z_a\subset \Delta^{n+k}\times X$. Evaluation map 
 (1.3.1) will be used in the proof of the Lemma 1.3.
 
 Recall that we suppose that our complex space $X$ is equipped with some 
  Hermitian metric $h$.

\smallskip\noindent\bf
Lemma 1.3. {\it Let a holomorphic map $f:\bar\Delta^n\times \bar A^k(r,1)\to X$ into 
a complex space $X$ be given. Suppose that:

1) for every $z\in \bar\Delta^n$ 
the restriction $f_z$ extends meromorphically onto the whole $k$-disk 
$\bar\Delta^k_z$;

2)  the volumes of graphs of these extensions are uniformly 
bounded;

3) there exists a compact $K\Subset X$ which contains $f(\bar\Delta^n\times \bar A^k(r,1))$ and 
$f(\bar\Delta_z^k)$ for all $z\in \bar\Delta^n$.

\noindent Then $f$ extends meromorphically onto $\Delta^{n+k}$. 
}

 \smallskip\noindent\bf
 Proof. \rm Denote by $\nu =\nu (K)$ the minimal volume of a compact $k$-dimensional analytic
 subset in $K$, $\nu >0$ by Lemma 2.3.1 from [Iv-4]. Denote by $W$ the maximal open
 subset of $\Delta^n$ such that $f$ extends meromorphically  onto
 $\Delta^n\times A^k(r,1)\cup W\times \Delta^k$. Set $S=\Delta^n\setminus
 W$. Let
 $$
 S_l=\{ z\in S: {\sl vol}(\Gamma_{f_z})\le l\cdot {\nu \over 2}\}.\eqno(1.3.2)
 $$
 The maximality of $W$ (and thus the minimality of $S$) and  Lemma 2.4.1 from [Iv-4] 
 imply that $S_{l+1}\setminus S_l$ are pluripolar
 and by the  Josefson theorem so is $S$. In particular, $W\not= \emptyset $.

 Consider the analytic space

 $$
 {\cal G}_{f,2C_0,c}:=\{ Z\in {\cal G}_{f,2C_0}: \Vert \phi (Z)\Vert <c\}
 ,\eqno(1.3.3)
 $$
 \noindent
 where $0<c\le 1$ is fixed. $C_0$ is taken here such that ${\sl vol}(\Gamma_{f_z})\le 
 C_0$ for all $z\in \bar\Delta^n$. Since, by Lemma 1.2 cycles of the form $\Gamma_{f_z}$ are 
 dense in ${\cal G}_{f,2C_0,1}$, we have that for every $Z\in {\cal G}_{f,2C_0,1}$ 
 in fact ${\sl vol}({\bf ev}(Z))\le C_0$. Therefore we see that 
 $\bar {\cal G}_{f,C_0,1}\cap  \phi^{-1}(\Delta^n(1))$
 is closed and open in ${\cal G}_{f,2C_0,1}$ and in fact coincides with ${\cal G}_{f,2C_0,1}$.
 Closures we take in the cycle space ${\cal G}_{f}$.
 
 For any $c<1$ the set $\bar{\cal G}_{f,C_0,c}=\phi^{-1}(\bar\Delta^n(c))$ is compact by the 
 Harvey-Shiffman generalization of Bishop's theorem. Therefore $\phi :{\cal G}_{f,2C_0,1}\to 
 \Delta^n$ is proper. Therefore ${\bf ev}:{\cal Z}_f\to \Delta^{n+k}\times X$ is also proper
 and by the Remmert proper mapping theorem its image is an analytic set extending the 
 graph of $f$. The latter follows from the fact that $\phi ({\cal G}_{f,2C_0,1})\supset W$
 and therefore in fact $\phi ({\cal G}_{f,2C_0,1})=\Delta^n(1)$.

\smallskip\hfill{q.e.d.}

 \smallskip\noindent\bf
 Definition 1.3. {\it We say that a complex space $X$ is disk-convex  in 
 dimension $k$ if for every compact $K\Subset X$  there exists a compact $\hat K$ 
 such that  for every meromorphic mapping $\phi :\bar\Delta^k\to X$ with 
 $\phi (\partial \Delta^k )\subset K$ one has $\phi (\bar\Delta^k)\subset \hat K$.
 }
 
 \smallskip\noindent\bf
 Remark 1. \rm For $k=1$ we say simply that $X$ is disk-convex.
 
 \noindent\bf
 2. \rm Recall that a complex space $X$ is called $k$-convex (in the sence 
 of Grauert) if there is an exhaustion function $\phi :X\to [0,+\infty [$ 
 which is $k$-convex at all points outside some compact $K$, i.e., its Levi 
 form has at least ${\sl dim}X-k+1$ positive eigenvalues. By an appropriate 
 version of the maximum principle for $k$-convex functions  
 $k$-convexity implies disk-convexity in dimension $k$.

 \noindent\bf
 3. \rm Condition (3) of  Lemma 1.3 (as well as of Theorems 1.4 and 1.5 below) is 
automatically satisfied if $X$ is supposed to be disk-convex in dimension $k$.

\smallskip\noindent\sl
1.4. The Levi-type extension theorem.

  \smallskip\rm
  In the proof of the Main Theorem we will deal with the situation where a holomorphic 
  map $f:\Delta^n\times A^k(r,1)\to X$ extends from  $A^k_z(r,1)$ to $\Delta^k_z$ not for all 
  $z\in \Delta^n$ but only for $z$ in some "thick" set $S$.
  
  \smallskip\noindent\bf Definition 1.4. 
  {\it A subset $S\subset \Delta^n$ we call thick at origin if for any 
  neighborhood $U$ of zero $U\cap S$ is not contained in a proper anlytic subset of $U$.
    
  \smallskip\rm Case of dimension two, i.e., $n=1$ is somewhat special. Let us consider this case 
  separately. In this case $S$ is thick at origin iff $S$ contains a sequence $\{ s_n\} $ which 
  converge to zero. 
  \smallskip\noindent\bf  Theorem 1.4. {\it Let $f:\Delta \times A(r,1)\to X$ be a holomorphic 
  map into a normal, 
  reduced complex space $X$. Suppose that for a sequence $\{ s_n\} $ of points in $\Delta $, 
  converging to the origin the restrictions $f_{s_n}:=f\mid_{A_{s_n}}$ extend holomorphically 
  onto $\Delta_{s_n}$. Suppose in addition that:
  
  1) there exists a compact $K\Subset X$ such that $\bigl[\bigcup_{n=1}^{\infty }f(\Delta_{s_n})\bigr]\cup 
  f(\Delta \times A(r,1))\subset K$;
  
  2) areas of images $f(\Delta_{s_n})$ are uniformly bounded.
  
  \noindent Then there exists an $\eps >0$ such that $f$ extends as a meromorphic map onto $\Delta (\eps )
  \times \Delta $.
  }
   
 \smallskip\rm  In dimensions bigger then two the situation becames more complicated, see examples 
 3.2 and 3.3 in \S 3. Let us give a condition on $X$ sufficient to maintain the conclusion of 
 Theorem 1.4. Denote by ${\bf ev}:{\cal Z}\to X$ the natural evaluation map from the 
 universal space ${\cal Z}$ over ${\cal B}_k(X)$ to $X$. 
 
 \smallskip\noindent\bf
 Definition 1.5. {\it Let us say that $X$ has unbounded cycle geometry in 
 dimension $k$ if there exists a path $\gamma :[0,1[ \to {\cal B}_k(X)$ with 
 ${\sl vol}_{2k}({\bf ev}(Z_{\gamma (t)}))\to \infty $ as $t\to \infty $ and ${\bf ev}(Z_{
 \gamma (t)})\subset K$ for all $t$, where $K$ is some compact in $X$.
 }

\smallskip\rm Now we can state the following

 \smallskip\noindent\bf
 Theorem 1.5. {\it Let $f:\Delta^n\times A^k(r,1)\to X$ be a holomorphic 
 mapping into a normal, reduced complex space $X$. Suppose 
 that there is a constant $C_0<\infty $ and a compact $K\Subset X$ such that 
 for $s$ in some subset $S\subset \Delta^n$, which is thick at origin the following holds:
 
 (a) the restrictions $f_s:=f\mid_{A^k_s(r,1)}$ extend meromorphically onto the 
 polydisk $\Delta^k_s$, and ${\sl vol}(\Gamma_{f_s})\le C_0$ for all $s\in S$;
 
 (b) $f(\Delta^n\times A^k(r,1))\subset K$ and $f_s(\Delta^k)\subset K$ for 
 all $s\in S$.
 
 \noindent
  If $X$ has bounded cycle geometry in dimension $k$, then 
 there exists a neighborhood $U\ni 0$ in $\Delta^n$ and a meromorphic extension
 of $f$ onto $U\times \Delta^k$. 
 }     
 \smallskip\rm We shall use the Theorem 1.5 in the case when $k=1$. In this case it admits a 
 nice refinement. A $1$-cycle $Z=\Sigma_jn_jZ_j$ is called rational if all $Z_j$ are rational 
 curves, i.e., an images of the Riemann sphere $\cc\pp^1$ in $X$ under a non constant holomorphic 
 mappings. Considering the space of rational cycles ${\cal R}(X)$ instead of 
 Barlet space ${\cal B}_1(X)$ we can define as in Definition 1.5 the notion of bounded 
 rational cycle geometry.
 
 \smallskip\noindent\bf Corollary 1.6. \it Suppose that in the conditions of Theorem 1.5  one has 
 additionally that 
 $k=1$. Then the conclusion of this theorem holds provided $X$ has bounded rational 
 cycle geometry.
 
 \smallskip\noindent\bf
 Proof of Theorems 1.4, 1.5 and Corollary 1.6.  \rm
 \smallskip\noindent\sl
 Case $n=1$. \rm  Define ${\cal G}_0$  as the set of all limits $\{ \Gamma_{f_{s_n}},s_n\in S,s_n\to 0 
 \} $. Consider the union $\hat {\cal G}_0$ of those components
 of ${\cal G}_{f,2C_0}$ that intersect ${\cal G}_0$. At least one of these 
 components, say ${\cal K}$, contains two points $a_1$ and $a_2$ such that 
 $Z_{a_1}$ projects onto $\Delta_0^k$ and $Z_{a_2}$ projects onto
 $\Delta_s^k$ with $s\not= 0$. This is so because $S$ contains a sequence 
 converging to zero. Consider the restriction ${\cal Z}_f\mid_{{\cal K}}$ of the
 universal family onto ${\cal K}$. This is a
 complex space of  finite dimension. Join the points $a_1$ and $a_2$ by an analytic 
 disk $h :\Delta \to  {\cal K}$, $h (0)=a_1, h (1/2)=a_2$. Then the
 composition $\psi = \phi\circ h :\Delta \to \Delta $ is not
 degenerate because $\psi (0)=0\not= s=\psi (1/2)$. 
 Here $\phi :=p_1\circ {\bf ev}\circ \pi^{-1}:{\cal C}_f\to \Delta^n$ was defined in 
 the proof of Lemma 1.2. Map $\phi $ restricted to ${\cal G}_f$ will be denoted also as 
 $\phi $. 
 Thus $\psi $ is proper
 and obviously so is the map ${\bf ev}:{\cal Z}\mid_{\psi (\Delta )}\to 
 F({\cal Z}\mid_{\psi (\Delta )})\subset \Delta^{1+k}\times
 X$. Therefore  ${\bf ev}({\cal Z}\mid_{\psi (\Delta )})$ is an analytic set in 
 $U\times  \Delta^k\times X$ for small enough $U$ extending $\Gamma_f$ by the 
 reason
 of dimension.
 
 \smallskip This proves the   {\sl Theorem 1.4}.
 \smallskip\noindent\sl
 Case $n\ge 2$. \rm We shall treat this case in two steps.

 \smallskip\noindent\sl
 Step 1. \it Fix a point $z\in \Delta^n$ such that 
 $\phi ({\cal G}_f)\ni z$. Then there exists a relatively compact open  $W\subset {\cal G}_{f}$,
 which contains ${\cal G}_{f,C_0}$ 
 such that $\phi (W)$ is an analytic variety in some neighborhood $V$ of $z$. 
 
 \smallskip\noindent\rm Consider the analytic subset $\phi^{-1}(z)$ in ${\cal G}_f$. Every 
 $Z_a$ with $a\in \phi^{-1}(z)$ has the form $B_a+\Gamma_{f_z}$, where $B$ is 
 a compact cycle in $\Delta^k_z\times X$. Thus connected components of 
 $\phi^{-1}(z)$ parametrize connected and closed subvarietes in ${\cal B}_k(
 \Delta^k\times X)$. Holomorphicity of $f$ on $\Delta^n\times A^k(r,1)$ and 
 condition (b) of the Theorem 1.5 imply that $B_a\subset \bar\Delta^k_z\times K$. So, 
 if $\phi^{-1}(z)$ had non compact connected components, this would imply 
 the unboundness of cycle geometry of $X$. 
 
 Thus, all connected components of  $\phi^{-1}(z)$  should be compact. Let ${\cal 
 K}$ denote the union of connected components of $\phi^{-1}(z)$ intersecting
 ${\cal G}_{f,C_0}$. Since ${\cal K}$ is compact, there obviously exist a 
 relatively compact open $W\Subset {\cal G}_f$ containing  ${\cal G}_{f,C_0}$
 and ${\cal K}$,  
 and a neighborhood $V\ni z$ such that $\phi\mid_W:W\to V$ is proper. 
 By the Remmert's proper mapping theorem. 
 $\phi (W)\subset V$ is an analytic subset of $V$. 
 
 \smallskip\noindent\sl
 Step 2. \it If $S$ is thick at $z$ then there exists a neighborhood $V\ni z$ 
 such that $f$ meromorphically extends onto $V\times \Delta^k$.
 
 \smallskip\noindent\rm  
 Since $\phi (W)\supset S\cap V$ and $S$ is thick at the origin, the first step implies that 
 $\phi (W)\cap V=V$. Since $W\Subset {\cal G}_f$ there exist a constant $C$ s.t. 
${\sl vol} \{ Z_s:s\in W\} \le 
 C$. This allows to apply the Lemma 1.3 and obtain the extension of $f$ onto $V\times \Delta^k$.

 \smallskip This proves the Theorem 1.5.
 \smallskip\noindent\sl
 Case $k=1$. \rm The limit of a sequence of analytic disks of bounded area is an 
 analytic disk plus a rational cycle, see ex. [Iv-1]. Therefore we need to 
 consider only the space of ratioanl cycles in this case. The rest is obvious. This 
gives Corollary 1.6. 
\smallskip
\hfill{q.e.d.}
 
 \smallskip Step 1 in the proof of Theorem 1.5 in fact gives the following statement, 
which will be used later.

 \smallskip\noindent\bf
 Corollary 1.7. {\it Let $f:\Delta^n\times A^k(r,1)\to X$ be a holomorphic 
 mapping into a normal, reduced complex space $X$ which has bounded cycle 
 geometry in dimension $k$. Suppose 
 that there is a constant $C_0<\infty $ and a compact $K\Subset X$ such that 
 for $s$ in some subset $S\subset \Delta^n$, the following holds:
 
 (a) the restrictions $f_s:=f\mid_{A^k_s(r,1)}$ extend meromorphically onto the 
 polydisk $\Delta^k_s$, and ${\sl vol}(\Gamma_{f_s})\le C_0$ for all $s\in S$;
 
 (b) $f(\Delta^n\times A^k(r,1))\subset K$ and $f_s(\Delta^k)\subset K$ for 
 all $s\in S$.
 
   Then there exists a neighborhood $V\ni 0$ and an analytic subvariety $W$ of 
 $V$ such that $W\supset S\cap V$ and such that for every $z\in W$ $f_z$ 
  meromorphically extends onto $\Delta^k_z$ with ${\sl vol}(\Gamma_{f_z})\le C_0$.
}




\smallskip

\rm In the same spirit one obtains the following
  
\smallskip\noindent\bf
Corollary 1.8. \it Let a meromorphic mapping $f:\Delta^n\times A(r,1)\to X$
be given, where $X$ is a compact complex manifold with bounded rational cycle 
geometry. Let $S$ be a subset of $\Delta^n$ consisting of such points
$s$ that $f_s$ is well defined and extends holomorphically onto $\Delta_s$.
If $S$ is not contained in a countable union of locally closed proper
analytic subvarietes of $\Delta^n$, then there exists an open non-empty
$U\subset \Delta^n$ and a meromorphic extension of $f$ onto $U\times \Delta $.

  \smallskip\rm

Indeed, one easily deduces the  existence of a point $p\in \Delta^n $, that
can play the role of the origin in Theorem 1.5.

\bigskip\noindent\sl
 1.5. A remark about spaces with bounded cycle geometry.
 \smallskip\rm To apply the Theorem 1.5 in the proof of the Main Theorem we need to 
 check the boundedness of cycle geometry of the manifold $X$ which carries a pluriclosed 
 metric form. We shall do this in Proposition 1.9 below.
 We start from the following simple observation: 
 
 \smallskip\noindent\it
 Every compact complex manifold of dimension $k+1$ carries a strictly positive 
 $(k,k)$-form $\Omega^k$ with $dd^c\Omega^k=0$.
 
 \smallskip\rm Indeed: either a compact complex manifold  carries a 
 $dd^c$-closed strictly positive  $(k,k)$-form or it  carryies a 
 bidimension $(k+1,k+1)$-current $T$ with $dd^cT\ge 0$ but $\not\equiv 0$. 
 In the case of $\dim X=k+1$ such current is nothing but a nonconstant 
 plurisubharmonic  function, which doesn't exists on compact $X$.

Let us introduce the class ${\cal G}_k$ of normal complex spaces,
carrying
a nondegenerate positive $dd^c$-closed strictly positive $(k,k)$-forms.
Note that the sequence
$\{ {\cal G}_k\}$ is rather exaustive: ${\cal G}_k$ contains all compact
complex
manifolds of dimension $k+1$.
 
 Introduce furthermore the class of normal complex spaces ${\cal P}_k^{-}$ which 
 carry a strictly positive $(k,k)$-form $\Omega^{k,k}$ with $dd^c\Omega^{k,k}
 \le 0$. Note that ${\cal P}_k^{-}\supset {\cal G}_k$. As was mentioned in the 
 Introduction a Hopf three-fold $X^3=\cc^3\setminus \{ 0\} /(z\sim 2z)$ belongs to 
  ${\cal P}_1^{-}$ but not to ${\cal G}_1$.

 \smallskip\noindent\bf
 Proposition 1.9. {\it Let $X\in {\cal P}_k^{-}$ 
 and let ${\cal K}$ be an irreducible component of ${\cal B}_k(X)$ such 
 that ${\bf ev}({\cal Z}\mid_{\cal K}) $ is relatively compact in $X$. Then:
 \smallskip
 1) ${\cal K}$ is compact.
 \smallskip
 2) If $\Omega^{k,k}$ is a $dd^c$-negative $(k,k)$-form on $X$, then
 $\int_{Z_s}\Omega^{k,k}\equiv const$ for $s\in {\cal K}$.
 \smallskip
 3) $X$ has bounded cycle geometry in dimension $k$.
 \smallskip\noindent\sl
 Proof. \rm Let ${\bf ev}:{\cal Z}\mid_{{\cal K}}\to X$ be the evaluation map, and
 let $\Omega^{k,k}$ be a strictly positive $dd^c$-negative $(k,k)$-form on
 $X$. Then $\int_{Z_s}\Omega^{k,k}$ measures the volume of $Z_s$. Let us
 prove that the function $v(s)=\int_{Z_s}\Omega^{k,k}$ is plurisuperharmonic
 on ${\cal K}$. Take an analytic disk $\phi :\Delta \to {\cal K}$. Then 
 for any nonnegative test function  $\psi $ on $\Delta $    by
 Stokes theorem and reasons of bidegree we have
 $$
 <\psi ,\Delta \phi^*(v)> = \int_{\Delta }\Delta \psi \cdot \int_{
 Z_{\phi (s)}}\Omega^{k,k} = \int_{{\cal Z}\mid_{\phi (\Delta )}}dd^c(\pi^*
 \psi )\wedge \Omega^{k,k} =
 $$
 $$
 = \int_{{\cal Z}\mid_{\phi (\Delta )}}\pi^*\psi \wedge dd^c\Omega^{k,k}
 \le 0.
 $$
 \noindent
 Here $\pi :{\cal Z}\mid_{{\cal K}}\to {\cal K}$ is the natural projection.
 So $\Delta \phi^*(v)\le 0$ for any analytic disk in ${\cal K}$ in the sence
 of distributions. Therefore $v$ is plurisuperharmonic.

 Note that by Harvey-Shiffman generalization of Bishop's theorem $v(s)\to
 \infty $ as $s\to \partial {\cal K}$. So by the minimum principle
 $v\equiv const$ and ${\cal K}$ is compact again by Bishop's theorem.
 \smallskip\noindent
 2) The same computation shows that $\int_{Z_s}\Omega^{k,k}$ is
 plurisuperharmonic for any $dd^c$-negative $(k,k)$-form. Since ${\cal K}$
 is proved to be compact, we obtain the statement.
 \smallskip\noindent
 3) Let ${\cal R}$ be any connected component of ${\cal B}_k(X)$. Write
 ${\cal R}=\bigcup_j {\cal K}_j$, where ${\cal K}_j$ are irreducible 
components. From (1) we have that $v$ is constant on
 ${\cal R}$. So if $\{ {\cal K}_j\} $ is not finite then ${\cal R}$ has an
 accumulation point $s=\lim s_j$ by Bishop's theorem, where all $s_j$ belong to different components
 ${\cal K}_j$ of ${\cal R}$. This contradicts  the fact that ${\cal B}_k(X)$
 is a complex space.
 \smallskip
 \hfill{q.e.d.}
\smallskip

\bigskip\noindent\bf
2. Hartogs-type extension and spherical shells.
\smallskip\noindent\sl
2.1. Generalities on pluripotential theory. \rm

For the standard facts from  pluripotential theory we refer to [Kl].
Denote by ${\cal D}^{k,k}(\Omega )$ the space of $C^{\infty }$-forms of
bidegree $(k,k)$ with compact support on a complex manifold $\Omega $. $\phi \in
{\cal D}^{k,k}(\Omega )$ is real if $\bar \phi =\phi $.
The dual space ${\cal D}_{k,k}(\Omega )$ is the space of currents of
bidimension $(k,k)$ (bidegree $(n-k,n-k)$, $n={\sl dim }_{\cc }\Omega $).
$T\in {\cal D}_{k,k}(\Omega )$ is real if $<T,\bar \phi >=\overline{<T,\phi >}$
 for all $\phi \in {\cal D}^{k,k}(\Omega )$.
\smallskip\noindent\bf
Definition 2.1. {\it A current $T\in {\cal D}_{k,k}(\Omega )$ is called
positive if for all $\phi_1,...,\phi_k\in {\cal D}^{1,0}(\Omega )$
$$
<T,{i\over 2}\phi_1\wedge \bar\phi_1\wedge ...\wedge {i\over 2}\phi_k\wedge
\bar\phi_
k>\ge 0.
$$
\smallskip\noindent
$T$ is negative if $-T$ is positive.
}
\smallskip\noindent\bf
Definition 2.2. {\it We  say that a current $T\in {\cal D}_{k,k}
(\Omega )$ is pluripositive (-negative) if $T$ is positive and $dd^cT$ is
positive (-negative).
$T$ is pluridefinite if it is either pluripositive or plurinegative.
A current $T$ (not necessarily positive) is pluriclosed
 if $dd^cT=0$.
 }
 
\smallskip\rm 

If $K$ is a complete pluripolar compact in strictly pseudoconvex domain $\Omega
\subset \cc^n$ and $T$ is a   closed, positive current on $\Omega \setminus K$, then
$T$ has locally finite mass in a neighborhood of $K$, see [Iv-2], {\sl Lemma
 2.1 }. For a current $T$, which has locally finite mass in a neighborhood of 
 $K$, one denotes by $\tilde T$ its trivial extension onto $\Omega $, see [Lg].
 
\smallskip\noindent\bf
Lemma 2.1. {\it (a) Let $K$ be a complete pluripolar compact in a strictly
pseudoconvex domain $\Omega \subset \cc^n$ and $T$ be a pluridefinite current
of bidegree (1,1)  on $\Omega \setminus K$ of locally finite mass in a
neighborhood of $K$ and such that $dT$ has coefficients measures in $\Omega
\setminus K$. Then $dd^c\tilde T$ has coefficients measures on $\Omega $.

\noindent
(b) If $n=2$ and $K$ is of Hausdorff dimension zero, then 
$\chi_K\cdot dd^c\tilde T$ is negative, where $\chi_K$ is the characteristic 
function of $K$.
}

\smallskip\noindent\sl
Proof. \rm Part (a) of this lemma was proved in [Iv-2], Proposition 2.3 for
currents of bidimension (1,1) (the
condition on $dT$ was forgotten there). If $T$ is of bidegree (1,1), then
consider $T\wedge (dd^c\Vert z\Vert^2)^{n-2}$ to get the same conclusion.

\smallskip\noindent
(b) Let $\{ u_k\} $ be a sequence of smooth plurisubharmonic functions in 
$\Omega $,
equal to zero in a neighborhood of $K$, $0\le u_k\le 1$ and such that
$u_k\nearrow \chi_{\Omega \setminus K}$ uniformly on compacts in $\Omega
\setminus K$, see {\sl Lemma 1.2} from [Sb]. Put $v_k=u_k-1$.

$\widetilde{dd^cT}$ is a negative measure on $\Omega $ , which we shall denote as $\mu_0$. 
According to
part (a) the distribution  $\mu :=dd^c\tilde T$
is a measure. Write
$$
\mu =\chi_K\cdot \mu + \chi_{\Omega \setminus K}\cdot \mu ,\eqno(2.1.1)
$$
\noindent where obviously $\chi_{\Omega \setminus K}\cdot \mu =\mu_0$. Denote 
the measure $\chi_K\cdot \mu $  by $\mu_s$. We shall prove that the 
measure $\mu_s$ is
non-positive.
Take a ball $B$ in $\cc^2$ centered at $s_0\in K$ such that $\partial B\cap K=
\emptyset $. One has
$$
\mu_s(B\cap K) = -\lim_{k\rightarrow \infty }\int_Bv_k\cdot \mu  =
- \lim_{k\rightarrow \infty }<v_k,dd^c\tilde T> 
= -\lim_{k\rightarrow \infty }<dd^cv_k,\tilde T> \le 0, \eqno(2.1.2)
$$
\noindent
because $\tilde T$ is positive and $dd^cv_k\ge 0$. So for any such ball we
have
$$
\mu_s(B\cap K)\le 0. \eqno(2.1.3)
$$
All that is left, is to use the following Vitali-type theorem for  general
measures,
see [Fd], p.151. Let $D$ be an open set in $\cc^2$ and $\sigma $ a finite
positive Borel measure on $D$. Let further ${\cal B}$ be a family of closed
balls of positive radii such that for any point $x\in D$ the family ${\cal B}$
 contains  balls of arbitrarily small radii centered at $x$. Then one can
find a countable subfamily $\{ B_i\} $ of pairvise disjoint balls in ${\cal B}$ 
such that
$$
\sigma (D\setminus \bigcup_{(i)}B_i) = 0. \eqno(2.1.4)
$$
Represent our measure $\mu_s$ as a difference $\mu_s=\mu_s^+ - \mu_s^-$ of two
nonnegative measures. Fix a relatively compact open subset $D\subset \Omega
$. As ${\cal B}$ take the family of all balls such that $\partial B\cap K=\emptyset
$.
 Since  $K$ is of dimension zero this is a Vitali-type covering. Let $\{ B_i
\} $ be pairwise disjoint and such that $\mu_s^+ (D\setminus \bigcup_{(i)}B_i) = 0$. 
Then
$\mu_s^+(D) = \mu_s^+(D\setminus \bigcup_{(i)}B_i) + \sum_{(i)}\mu_s^+(B_i) =
\sum_{(i)}\mu_s^+(B_i)$. Consequently
$$
\mu_s(D) = \mu_s^+(D) - \mu_s^-(D) \le \mu_s^+(\bigcup_{(i)}B_i) -
\mu_s^-(\bigcup_{(i)}B_i) =
$$
$$
= \sum_i\mu_s^+(B_i) - \sum_i\mu_s^-(B_i) = \sum_i\mu_s(B_i)\le 0 \eqno(2.1.5)
$$
\noindent by (2.1.3). Thus $\mu_s(D)\le 0 $ for any relatively compact open
set $D$ in $\Omega $. So
the measure $\mu_s$ is negative.
\smallskip

\hfill{q.e.d.}

\smallskip\rm
Together with Main Theorem we shall prove a somewhat more general result. A 
metric form $\omega $ on $X$ we call plurinegative if $dd^c\omega \le 0$. But first 
recall the following

\smallskip\noindent\bf
Definition 2.3. {\it Recall that a subset $K\subset \Omega $ is called (complete)
$p$-polar if for any $a\in \Omega $ there exist a neighborhood $V\ni a$
and  coordinates $z_1,...,z_n$ in $V$ such that the sets
$
K_{z^0_I} = K\cap \{ z_{i_1}=z_{i_1^0},...,z_{i_p}=z_{i_p^0}\} 
$
are (complete) pluripolar in the subspaces $V_{z_i^0}:=\{ z\in V:
 z_{i_1}=z_{i_1^0},...,z_{i_p}=z_{i_p^0}\} $ for almost all $z_I^0=(z_{i_1}^0,
...,z_{i_p}^0)\in \pi^I(V)$, where $I$ runs over a finite set of multiindices
with $\vert I\vert =p$, such that $\{ (\pi^{I})^*w_e\} _I $ generates the
space of $(p,p)$-forms. Here $\pi^I(z_1,...,z_n)=(z_{i_1},...,
z_{i_p})$ denotes the projection onto
the space of variables $(z_{i_1},...,z_{i_p})$.
}

\smallskip\rm Now we can state the main result in most general form.

\smallskip\noindent\bf
 Theorem 2.2. {\it Let $f:H_U^{n+1}(r)\to X$ be a meromorphic map into a disk-convex
  complex space $X$ that admits a plurinegative Hermitian metric form $\omega $. Then:

  (1) $f$ extends to a meromorphic map  $\hat f:\Delta^{n+1}\setminus A
  \to X$, where $A$ is a closed, complete $(n-1)$-polar subset of $\Delta^{n+1}$ of Hausdorff 
  $(2n-1)$-dimensional measure zero.

  (2) If, in addition, $\omega $ is pluriclosed and if $A\not= \emptyset
   $ is the minimal subset such that $f$ extends onto $\Delta^{n+1}\setminus A$, then for 
   every transversal sphere $\ss^3\subset \Delta^{n+1}
  \setminus A$ its image $f(\ss^{3})$  is not homologous to zero
  in $X$.
}

\smallskip\rm We would like to turn attention to the difference between plurinegative and 
pluriclosed cases. Example of the Hopf three-fold, given in the Introduction, shows that 
when $X$ admits only 
plurinegative metric form the singular set $A$ can have "components" of Hausdorff codimension 
higher 
than four and that the homological characterization  of $A$ is also not valid in 
general.

\bigskip\noindent\sl
2.2. Proof in dimension two.

\smallskip\rm
Let a meromorphic mapping $f:H_U^2(1-r)\to X$ from the two-dimensional Hartogs
figure into a disk-convex complex space be given. Since the indeterminancy set $I(f)$ of $f$ is 
discrete, we can suppose after shrinking $A(1-r,1)$ and $\Delta $ if necessary, that $f$ is 
holomorphic in the neighborhood of $\bar\Delta \times \bar A(1-r,1)$. Let $\omega $ be a plurinegative
metric form on $X$. Denote by $W$  the maximal open subset of the unit disk
$\Delta $ such that $f$ extends holomorphically  onto $H_W^2(1-r):=W\times \Delta
\cup \Delta \times A(1-r,1)$. Note that $W$ contains $U$ exept possibly a discrete set. Let $I(f)$ be 
the fundamental set of $f$ and denote by $\hat f$  the mapping $\hat f(z)=(z,f(z))$ into the graph. 
For $z\in W$ define
$$
a(z) = {\sl area} \hat f(\Delta_z) = \int_{\Delta_z}(dd^c\vert
\lambda \vert ^2+f\mid_{\Delta_z}^*\omega ).\eqno(2.2.1)
$$
\noindent
Here $\Delta_z=\{ (z,\lambda ):\vert \lambda \vert <1\} $.
We start with the following simple observation. Denote by $\nu_1=\nu_1(K)$
 the infimum of areas of compact complex curves contained in a 
compact $K\Subset X$. Then $\nu_1>0$, see {\sl Lemma 2.3.1} in [Iv-4].
\smallskip\noindent\bf
Lemma 2.3. \it Let $f:\bar\Delta \times \bar A(1-r,1)\to X$ be a holomorphic mapping
into a disk-convex complex space $X$. Suppose that for some sequence of points
 $\{ s_n\} \subset \Delta $, $s_n\rightarrow 0$, the following holds:

(a) $f_{s_n}:=f\mid_{\{ s_n\} \times A(1-r,1)}$ extends holomorphically onto $\Delta
_{s_n}:=\{ s_n\} \times \Delta $;

(b) ${\sl area}\hat f(\Delta_{s_n})\le C$ for all $n$.

\noindent Then $f_0:=f|_{\{s_0\} \times A(1-r,1)}$ extends holomorphically 
onto $\Delta_0$.

\smallskip If moreover,

(c) for  a compact $K$  in $X$ containing the set $f[(\Delta ({1\over 2})\times A(1-
{2\over 3}\cdot r,1-{1\over 3}\cdot r)\cup \bigcup_{(n)}\{ s_n\} \times \Delta (1-{1\over 3}\cdot
r)]$ one has

$$
\bigl\vert {\sl area} \hat f(\Delta_{s_n}(1-{1\over 3}\cdot r)) - {\sl area} \hat
f(\Delta_0(1-{1\over 3}\cdot r)\bigr\vert \le {1\over 2}\cdot \nu_1(K), \eqno(2.2.2)
$$

\noindent
for $n\gg 1$, then $f$ extends holomorphically  onto $V\times \Delta $ for some open 
$V\ni 0$.
\smallskip\noindent\sl
Proof. \rm The first statement is standard. Let us prove the second one. 
 First of all let us show that ${\cal H}-\lim_{n\rightarrow \infty }
\hat f(\bar\Delta_{s_n}(1-{1\over 3}\cdot r))=\hat f(\bar\Delta_0(1-{1\over 3}\cdot r))$,
i.e., the sequence of graphs $\{ \hat f(\bar\Delta_{s_n}(1-{1\over 3}\cdot r))\} $
converges in the Hausdorff metric to the graph of the limit. If not, there would
be a subsequence (still denoted by $\{ \hat f(\bar\Delta_{s_n}(1-{1\over 3}\cdot r))
\} $) such that  ${\cal H}-\lim_{n\rightarrow \infty }
\hat f(\bar\Delta_{s_n}(1-{1\over 3}\cdot r))=\hat f(\bar\Delta_0(1-{1\over 3}\cdot r))\cup
\bigcup_{j=1}^N\{ p_j\} \times C_j$, where $\{ C_j\} $ are compact curves,
see {\sl Lemma 2.3.1} in [Iv-4]. Thus by (2.3.2) from [Iv-4] we have   
${\sl area} \hat f(\bar\Delta_
{s_n}(1-{1\over 3}\cdot r))\ge {\sl area }\hat f(\bar\Delta_0(1-{1\over 3}\cdot r))+N\cdot
\nu_1(K)$. This contradicts (2.2.2).

Take a Stein neighborhood $V$ of  $\hat f(\bar\Delta_0(1-{1\over 3}\cdot r))$, see
[Si-1]. Then for $\delta >0$ small enough we have  $f(\Delta_{\delta }
\times A_{1-{1\over 3}r-\delta ,1-{1\over 3}r+\delta })\subset V$ and $f(\Delta_{s_n}(1-{1\over 3}r)
)\subset V$ if $s_n\in \Delta_{\delta }$. From Hartogs theorem for holomorphic
functions we see that $f$ extends to a holomorphic map from $\Delta_{\delta }
\times \Delta_{1-{1\over 3}r-\delta }$ to $V$.
\smallskip
\hfill{q.e.d.}

\smallskip
\noindent
\bf Lemma 2.4. \it If the metric form $\omega $ on a disk-convex complex space $X$ 
is plurinegative and $W$ is maximal, then $\partial W\cap \Delta $ is complete polar 
in $\Delta $.
\smallskip\noindent\sl
Proof. \rm Take a point $z_0\in \partial W\cap \Delta $. Choose  a relatively
compact  neighborhood $V$ of $z_0$ in $\Delta $. Denote by $T={i\over 2}t^
{\alpha \bar\beta }dz_{\alpha }\wedge d\bar z_{\beta }$ the current $f^*\omega +dd^c
\Vert z\Vert^2$.
 The area function from (2.2.1) can be now written  as
$$
a(z_1) = {i\over 2}\cdot \int_{\vert z_2\vert \le 1}t^{2\bar 2}(z_1,z_2)dz_2\wedge
d\bar z_2 .\eqno(2.2.3)
$$
The condition that $dd^cT $ is negative means that
$$
{\partial ^2t^{1\bar 1} \over \partial z_2\partial\bar z_2}+
{\partial ^2t^{2\bar 2} \over \partial z_1\partial\bar z_1}-
{\partial ^2t^{1\bar 2} \over \partial z_2\partial\bar z_1}-
{\partial ^2t^{2\bar 1} \over \partial z_1\partial\bar z_2}\le 0 \eqno(2.2.4)
$$
\noindent
 on $H_W^2(1-r)$. Now we can estimate the Laplacian of $a $:
$$
\Delta a(z_1) = i\int_{\vert z_2\vert \le 1}
{\partial ^2t^{2\bar 2} \over \partial z_1\partial\bar z_1}dz_2\wedge d\bar z_2
\le i\int_{\vert z_2\vert \le 1}(-{\partial ^2t^{1\bar 1} \over \partial z_2
\partial\bar z_2}+{\partial ^2t^{1\bar 2} \over \partial z_2\partial\bar z_1} +
{\partial ^2t^{2\bar 1} \over \partial z_1\partial\bar z_2})dz_2\wedge d\bar z_2
 =
$$
$$
 = i\int_{\vert z_2\vert = 1}{\partial t^{1\bar 1} \over \partial z_2
}dz_2 + i\int_{\vert z_2\vert =1}{\partial t^{1\bar 2} \over \partial\bar z_1}
 d\bar z_2  - i\int_{\vert z_2\vert =1}{\partial t^{2\bar 1} \over \partial z_
1} dz_2 = \psi (z_1).\eqno(2.2.5)
$$
Inequality (2.2.5) holds for $z_1\in V\cap W$. But the right hand side $\psi $
is smooth in the whole $V$. Let $\Psi $ be a smooth solution of $\Delta \Psi
=\psi $ in $V$. Put $\hat a (z)=a(z)-\Psi (z)$. Then $\hat a $ is
superharmonic and bounded from below in $V\cap W$, maybe after shrinking $V$.

Denote further by $E$ the set of points $z_1\in \partial W\cap V$ such that
$a(z)\rightarrow +\infty $ as $z\in W,z\rightarrow z_1$. Note that $\hat
a(z)$ also tends to $+\infty $ in this case. For any point $z_{\infty }\in
[\partial W\cap V]\setminus E$ we can find a sequence $\{ z_n\} \subset W,z_n
\rightarrow z_{\infty }$ such that $a_t(z_n)\le C$. So by {\sl Lemma 2.3}
$f\mid_{\Delta_{z_{\infty }}\setminus \Delta_{z_{\infty }}(1-r)}$ extends onto
$\Delta_{z_{\infty }}$. 

Let $\nu_1$ be from {\sl Lemma 2.3} above for an appropriate $K\Subset X$. 
This compact $K$ should be taken to contain $f(\bar V\times \bar A(1-r,t)\cup
(W\cap \bar V)\times \bar\Delta_t)$. It exists because of disk-convexity of $X$. 
Set $E_j=\{ z\in \partial W\cap V:
a(z)\le {j\over 2}\nu_1\} $ for $j=1,2,...$. From {\sl Lemma 2.3} we
see that $E_j$ are closed subsets of $\partial W\cap V$, $E_j\subset E_{j+1}$
and we have that $\partial W\cap V = E\cup \bigcup_{j=1}^{\infty }E_j$.

Furthermore from {\sl Lemma 2.3}  we see that
$E_{j+1}\setminus E_j$ is a discrete subset of $V\setminus E_j$, say $E_{j+1}
\setminus E_j = \{ a_{ij}\} $. Now put
$$
u_1(z) =-\sum_{i,j}c_{ij}\log \vert z-a_{ji}\vert. \eqno(2.2.6)
$$
Here positive constants $c_{ij}$ are chosen in such a manner that
$\sum_{i,j}c_{ij}< +\infty $. Then $u_1(z)$ is superharmonic in $V$,
$u_1(z)\rightarrow +\infty $ as $z\rightarrow \bigcup_{j=1}^{\infty }E_j$
and $u_1(z)\not= +\infty $ for all $z\in V\cap W$. Now put $u_2(z)= \hat
a(z) + u_1(z)$. Note that $u_2$ is superharmonic in $W\cap V$ and $u_2(z)
\rightarrow +\infty $ as $z\rightarrow \partial W\cap V$. Define
$$
u_n(z) = \min \{ n, u_2(z)\} \eqno(2.2.7)
$$
\noindent
for $n\ge 3$. Note that $u_n$ are superharmonic in $V$, because $u_n\equiv n$
in
the neighborhood of $\partial W\cap V$. Put now $u(z) =\lim_{n\rightarrow
\infty }u_n(z)$. Then $u$ is superharmonic in $V$ as a nondecreasing limit of
superharmonic functions. Using the fact that $\hat a $ is
finite on $W$, we obtain that $u(z)=u_2(z)\not= +\infty $ for any $z\in V\cap
W$ and $u\mid_{V\setminus W}\equiv +\infty $, i.e.  $\partial W\cap \Delta $ is
  {\it complete} polar in $\Delta $.
So the lemma  is proved.
\smallskip
\hfill{q.e.d.}

\smallskip
In what follows we shall use the fact that a closed set of zero harmonic
measure in the plane has zero Hausdorff dimension, see [Gl]. Put $S_1=\Delta
 \setminus W$, where $W$ is the maximal domain in $\Delta $ such that our map $f$
extends holomorphically onto $H_W^2(1-r)$. We have proved that $S_1$ is polar i.e. 
of harmonic
measure zero. In particular, $S_1$ is zero dimensional. For any $\delta >0$ we
can find $0<\delta_1<\delta $ such that $\partial \Delta_{1-\delta_1}\cap S_1
=\emptyset $. Now we can change coordinates $z_1,z_2$ and consider the
Hartogs figure $H = \{ (z_1,z_2)\in \Delta^2: 1-r<\vert z_2\vert <1, \vert z_1
\vert <1 $ or $\vert z_2\vert <1, 1-\delta_1-\eps <\vert z_1\vert <1-\delta_1+
\eps \} $, where $\eps $ is small enough. Applying {\sl Lemma 2.4} again  we 
extend $f$ onto $\Delta \times (\Delta \setminus S_2)$ where $S_2$
is of harmonic measure zero. Therefore we obtain a holomorphic extension of $f$ 
onto $\Delta^2\setminus S$, where $S$ is a product of two complete polar sets in 
$\Delta $. So $S$ is complete polar itself and has Hausdorff dimension zero.
This proves Part 1 of the Theorem 2.2 in dimension two.

\smallskip 
Denote by $T$ the positive (1,1)-current (in fact the smooth form)
$f^*\omega $ on $\Delta^2\setminus S$. By {\sl Lemma 3.3} from [Iv-2] we have that
$T$ has locally summable coefficients on the whole $\Delta^2$ and from {\sl 
Lemma 2.1} above we see that $dd^c\tilde T$ is a negative measure with singular 
support contained in $S$.
We write $dd^c\tilde T=\mu $. We set 
furthermore $\mu_s:=\chi_S\cdot \mu $ and $\widetilde{dd^cT}=\mu_0$. 
All $\mu ,\mu_s $ and $\mu_0$ are negative measures, in fact $\mu_0$ is an 
$L^1$-function and $\mu = \mu_0+\mu_s$.

\smallskip Let us suppose now that the metric form $\omega $ on $X$ is 
pluriclosed. Shrinking, if necessary we shall suppose that $S$ is compact.

\smallskip\noindent\bf
Lemma 2.5. \it Suppose that the metric form $\omega$ is pluriclosed and take a
ball $B\subset\subset \Delta^2$ such that $\partial B\cap S=\emptyset $.

(i) If $f(\partial B)$ is homologous to zero in $X$ then $dd^c\tilde T=0$ on
$B$.

(ii) If $dd^c\tilde T=0$ then $f$ extends meromorphically  onto $B$.
\smallskip
 \rm In [Iv-2],  {\sl Lemma 4.4 }, this statement was proved for the case
when $S\cap B=\{ 0\} $. One can easily check that the same proof goes through for the
case when $S\cap B$ is closed zero dimensional. In fact in Lemmas 2.8 and 2.9 we will 
prove this statement  "with parameters".

So, statement $(2)$ of {\sl Theorem 2.2} and thus {\sl Main Theorem} are proved in 
the case $n=1$, i.e. in dimension two.

\medskip
\noindent\sl
2.3. Proof in higher dimensions: plurinegative metrics.

\rm Let us turn to the proof of these theorems in higher dimensions.
First we suppose that the metric form $\omega $ is 
plurinegative.

\smallskip \rm Let $f:H^{n+1}_U(1-r)\to X$ be our map. 
It will be convenient to set $U=\Delta^n(r)$. 
\smallskip
\noindent\sl
Step 1. \it $f$ extends to a holomorphic map of $\bigcup_{z^{'}\in\Delta^{n-1}
_r\setminus R_1}(\Delta_{z^{'}}^2\setminus S_{z^{'}})$ into $X$, where $R_1$ is
contained in a locally finite union of locally closed proper
subvarieties of $\Delta^{n-1}_r$ and $S_{z^{'}}$ is zerodimensional and
pluripolar in $\Delta^2_{z^{'}}$.
\smallskip\noindent\sl
Proof of Step 1. \rm
For $z'=(z_1,...,z_{n-1})\in \Delta^{n-1}_r$ denote by $H^2_{z'}=H^2_{z'}(1-r)$  the
two-dimensional Hartogs domain $\{ z'\} \times H^2(1-r)$ in the bidisk 
$\Delta^2_{z'}=
\{ z'\} \times \Delta^2\in \cc^{n+1}$. Shrinking $H^{n+1}(1-r)$ if necessary,
we can suppose that $I(f)$ consists of finitely many irreducible components.
Denote by $R_1$ the set of $z'\in \Delta^{n-1}_r$ such that ${\sl dim }[H^2_
{z'}\cap I(f)]>0$. $R_1$ is clearly contained in  finite union of
locally closed proper analytic subsets of $\Delta_r^{n-1}$. For $z'\in \Delta^
{n-1}_r\setminus R_1$, by the results of \S 2.2 the map $f\mid_{H^2_{z'}}$ 
extends to a
holomorphic map $f_{z'}:\Delta^2_{z'}\setminus S_{z'}\to X$, where $S_{z'}$ is
a zerodimensional and complete pluripolar in $\Delta^2_{z'}$. Note also  that 
$S_{z'}\supset \Delta^2_{z'}\cap I(f)$.

Take a point $z'\in \Delta^{n-1}_r\setminus R_1$ and a point $z_n\in \Delta
\setminus \pi_n(S_{z'})$. Here $\pi_n:\{ z'\} \times \Delta \times\Delta \to
\{ z'\}
\times \Delta $ is the projection onto the variable $z_n$. Take a domain $U
\subset\subset \{ z'\} \times \Delta \times \{ 0\} $ that is biholomorphic
to the unit disk, doesn't contain points from $\pi_n(S_{z'})$ and contain
the points $u:=(z',0,0) $ and $v:=(z',z_n,0)$. We also take  $U$ intersecting 
$A(1-r,1)$.
If $\{ z'\} \times \{ 0\}$ is
in $\pi_n(S_{z'})$ then take as $u$ some point close to $(z',0,0)$ in $\{ z'\}
\times \Delta $. Find a Stein neighborhood $V$ of the graph $\Gamma_{f\mid_
{\{ z'\} \times \bar U\times \Delta }}$. Let $w\in \partial U\cap A(1-r,1)$ be
some point. We have $f(\{ z',w\} \times \Delta )\subset V$ and
$f(\{ z'\} \times \partial U\times \Delta )\subset V$. So the usual
continuity principle for holomorphic functions
gives us a {\it holomorphic } extension of $f$ to the neighborhood of $\{
z'\} \times \bar
U\times \Delta $ in $\Delta^{n+1}$. Changing little bit the slope  of the 
$z_{n+1}$-axis 
and repeating the arguments as above  we obtain a holomorphic   extension of
$f$ onto the neighborhood of $\{ z'\} \times (\Delta \setminus S_{z'})$ for
each $z'\in \Delta^{n-1}_r\setminus R_1$.
\smallskip\noindent\sl
Step 2. \it $f$ extends holomorphically onto $(\Delta^{n-1}_r\times \Delta^2)
\setminus R$, where $R$ is a closed subset of $\Delta^{n-1}_r\times \Delta^2$
of Hausdorff codimension $4$.

\smallskip\noindent\sl
Proof of Step 2. \rm Consider a subset $R_2\subset R_1$ consisting of such
$z^{'}\in \Delta^{n-1}_r$ that $\dim [H^2_{z^{'}}\cap I(f)]=2$, i.e.
$H^2_{z^{'}}\subset I(f)$. This is a finite union of locally closed
subvarieties  of $\Delta^{n-1}_r$ of complex codimension at least two. Thus
$\bigcup_{z^{'}\in R_2}\Delta^2_{z^{'}}$ has Hausdorff codimension at least
four.

For $z^{'}\in R_1\setminus R_2=\{ z^{'}\in \Delta^{n-1}_r:\dim [H^2_{z^{'}}
(1-r)\cap I(f)]=1\} $ using \S 2.2 we can extend $f_{z^{'}}$
 holomorphically onto
$\Delta^2_{z^{'}}$ minus a zero dimensional polar set. Repeating the arguments
from {\sl Step 1} we can extend $f$ holomorphically  to a  neighborhood of
$\Delta^2_{z^{'}}\setminus C_{z^{'}}$ in $\Delta^{n-1}_r\times \Delta^2$.
Here $C_{z^{'}}$ is a complex curve containing all one dimensional components
of $H^2_{z^{'}}(1-r)\cap I(f)$.

$\bigcup_{z^{'}\in R_1\setminus R_2}C_{z^{'}}$ has Hausdorff codimension
at least four. Thus the proof of {\sl Step 2} is  completed by setting 
$R=\bigcup_{z^{'}\in R_1\setminus R_2}C_{z^{'}}\cup \bigcup_{z^{'}\in
R_2}\Delta^2_{z^{'}}$.

\smallskip\noindent\sl
Step 3. \rm We shall state this step in the form of the {\sl Lemma}.

\smallskip\noindent\bf
Lemma 2.6. \it There exists a closed, complete (n-1)-polar subset $A\subset R$ and a
holomorphic extension of $f$ onto $(\Delta^{n-1}_r\times \Delta^2)\setminus
A$ such that the current $T:=f^*\omega $ has locally summable coefficients in a
neighborhood of $A$. Moreover, $dd^c\tilde T$ is negative, where $\tilde T$ 
is the trivial extension of $T$.

\smallskip\rm
Take a point $z_0\in R$ and using the fact that $R$ is of Hausdorff
codimension four in $\cc^{n+1}$, find a neighborhood $V\ni z_0$ with a
coordinate system $(z_1,...,z_{n+1})$ such that $V=\Delta^{n-1}\times \Delta^2
$ in these coordinates and for all $z'\in \Delta^{n-1}$ one has  $R\cap
\partial \Delta^2_{z'}=0$. By \S 2.2
the restrictions $f_{z'}$ extend holomorphically onto $\Delta^2_{z'}\setminus
A(z')$, where $A(z')$ are closed complete pluripolar subsets in 
$\Delta^2_{z'}$ of
Hausdorff dimension zero. By the arguments similar to those used in  {\sl
Step 1 } $f$ extends holomorphically to a neighborhood of $V\setminus A$,
$A:=\bigcup_{z'\in \Delta^{n-1}}A(z')$.

Consider now the current $T=f^*\omega $ defined on $(\Delta^{n-1}\times \Delta^2)
\setminus R$. Note that $T$ is smooth, positive and $dd^cT\le 0$ there. By 
{\sl Lemma 3.3} from [Iv-2] every restriction $T_{z'}:=T\mid_{\Delta^2_
{z'}}\in L^1_{loc}(\Delta^2_{z'})$, $z'\in \Delta^{n-1}$. We shall use the
following Oka-type inequality for plurinegative currents proved in [F-Sb]:

\it there is a constant $C_{\rho }$ such that for any plurinegative current
$T$ in $\Delta^2$ one has
$$
\Vert T\Vert (\Delta^2) + \Vert dd^cT\Vert (\Delta^2)\le C_{\rho }\Vert T
\Vert (\Delta^2\setminus \bar\Delta^2_{\rho }).\eqno(2.3.1)
$$
\noindent
Here $0<\rho <1$. \rm

Apply (2.3.1) to the the trivial extensions $\tilde T_{z'}$ of $T_{z'}$, which are
plurinegative by ($b$) of {\sl Lemma 2.1}, to obtain that the masses $\Vert
\tilde T_{z'}\Vert (\Delta^2)$ are uniformly bounded on $z'$ on compacts in
$\Delta^{n-1}$. On $L^1$ the mass norm coincides with the $L^1$-norm. So
taking the second factor in $\Delta^{n-1}\times \Delta^2$ with different slopes
 and using Fubinis theorem we obtain that $T\in L^1_{loc}(\Delta^{n-1}\times
\Delta^2)$.

All that is left to prove is that $dd^c\tilde T$ is negative. It is enough to
show that for any collection $L$ of $(n-1)$ linear functions $\{ l_1,...,l_
{n-1}
\} $ the measure $dd^c\tilde T\wedge {i\over 2}\partial l_1\wedge
\overline{\partial l_1}\wedge ...\wedge {i\over 2}\partial l_{n-1}\wedge 
\overline{\partial l_{n-1}} $ is nonpositive, see [Hm]. Complete these functions to a coordinate
system
$\{ z_1=l_1,...,z_{n-1}=l_{n-1},z_n,z_{n+1}\} $ and note that for almost all
collections $L$  the set $\Delta^2_{z'}\cap A$ is of Hausdorff dimension
zero for all $z'\in \Delta^{n-1}$. Thus $\tilde T\mid_{z'}$ is plurinegative for all such $z'$. Take a
nonnegative function $\phi \in {\cal D}(\Delta^{n+1})$. We have 
$$
(n-1)!<dd^c\tilde T\wedge {i\over 2}\partial l_1\wedge
\overline{\partial l_1}\wedge ...\wedge {i\over 2}\partial l_{n-1}\wedge 
\overline{\partial l_{n-1}} ,\phi > = \int_{\Delta^{n+1}}
\tilde T\wedge (dd^c\Vert z'\Vert^2)^{n-1}\wedge dd^c\phi =
$$
$$
=  \int_{\Delta^{n-1}} (dd^c\Vert z'\Vert^2)^{n-1}
\int_{\Delta^2}(\tilde T)_{z'}\wedge  dd^c\phi =
\int_{\Delta^{n-1}} (dd^c\Vert z'\Vert^2)^{n-1}
\int_{\Delta^2}\tilde T_{z'}\wedge  dd^c\phi =
$$
$$
= \int_{\Delta^{n-1}} (dd^c\Vert z'\Vert^2)^{n-1}
\int_{\Delta^2}dd^c(\tilde T)_{z'}\wedge  \phi \le 0.
$$
We used here Fubini's theorem for $L^1$-functions, the fact that $(\tilde T)
_{z'}=\tilde T_{z'}$ for currents from $L^1_{loc}$ that are smooth outside
of a suitably situated set $A$, and finally the plurinegativity of $\tilde T
_{z'}$.

Therefore $\tilde T$ is plurinegative. We got an extension of $f$ onto $\Delta^{n-1}_r\times \Delta^2
\setminus A$, but this 
obviously  implies an extension (holomorphic) onto $\Delta^{n+1}
\setminus A$, where $A$ has zero Hausdorff  $(2n-2)$-dimensional measure 
 and is complete $(n-1)$-polar. 
 
 \smallskip\hfill{q.e.d.}

Therefore the 
first part of the Theorem 2.2 is proved.

The following statement is interesting by itself, but will be not used later. Let $W$ 
be the maximal open subset of $\Delta^{n+1}$ such that  $f$ meromorphically extend 
onto $W$. Denote by $\theta (x_0,dd^c\tilde T)$. the Lelong number of the closed negative
current $dd^c\tilde T$.

\smallskip\noindent\bf
Lemma 2.7. \it Under the conditions as above if  $x_0\in W$ 
then  $\theta (x_0,dd^c\tilde T)=0$.
\smallskip\noindent\bf Proof. \rm Find an orthonormal coordinate system $(z_1,...,z_{n-1})=z',(z_n,z_{n+1})=z''$ with center in 
$x_0$ and $r_0>0$ such that for every $x'\in \Delta^{n-1}(r_0)$ the intersection $\Delta^2_{x'}\cap {\cal I}(f) $ 
is finite. Here   $ {\cal I}(f) $ is the indeterminancy set of $f$. 
 For $r_0>r>0$ set $\tilde\Gamma_f(r):=p^{-1}(\Delta^{n+1}(x_0,r))$ and $\tilde\Gamma_{f,x'}(r):=
p^{-1}(\Delta^2_{x'}(r))$.
Note that by the geometric flattening theorem, see [Ba-2], 
$vol [\tilde\Gamma_{f,x'}(r_0)]\le C$ for all $x'\in \Delta^2(r_0)$. Now, because $\tilde T=T$ on $W$ we see 
that  

$$\theta (x_0,dd^c\tilde T)=
\lim_{r\to 0}{1\over r^{2(n-1)}}\int_{\Delta^{n+1}(x_0,r)} dd^cT\wedge (dd^c\Vert z\Vert^2 )^{(n-1)}.\eqno(2.3.6)
$$

\noindent This integral (and limit) can be estimated by the sum of integrals of the type 
$$
\lim_{r\to 0}{1\over r^{2(n-1)}}\int_{\Delta^{n+1}(x_0,r)} dd^cT\wedge (dd^c\Vert z'\Vert^2 )^{(n-1)},\eqno(2.3.7)
$$

\noindent taking sufficiently many orthonormal coordinate systems centered at $x_0$. So let us prove that the 
last limit is zero. First remark that 

$$
\int_{\Delta^{n+1}(0,r)} dd^cT\wedge (dd^c\Vert z'\Vert^2 )^{(n-1)}\le  
\int_{\Delta^{n-1}(0,r)}(dd^c\Vert z'\Vert^2)^{n-1}\int_{\Delta_{z'}^2(0,r)}dd^cT\le 
$$

$$
\le  r^{2(n-1)}\sup_{z'\in \Delta^{(n-1)}(0,r)}\int_{\Delta_{z'}^2(0,r)}dd^cT. \eqno(2.3.8)
$$
So we need to prove that $\sup_{z'\in \Delta^{(n-1)}(0,r)}\int_{\Delta_{z'}^2(0,r)}dd^cT\to 0$ as $r\to 0$. 
This will be done in two steps.

\smallskip\noindent\sl
Step 1. \it $\int_{\Delta_{0'}^2(0,r)}dd^cT\to 0$ as $r\to 0$.

\smallskip\rm Take the irreducible component $\Gamma_{f,0}(r_0)$ of $\tilde\Gamma_{f,0'}(r_0)$ which project onto 
$\Delta_{0'}^2(0,r)$ surjectively. This is the graph of the restriction $f|_{\Delta_{0'}^2(0,r_0)}$. Do the same 
for all $r<r_0$. Remark that  $\Gamma_{f,0}(r)\subset \Gamma_{f,0}(r_0)$. As $r\to 0$  $\Gamma_{f,0}(r)$ 
contracts to a finite union of curves - the fiber of $f|_{\Delta_{0'}(0,r_0)}$ over zero. In particular 
$vol[\Gamma_{f,0}(r)]\to 0$. Since $dd^cT=dd^c\pi^*\omega $ is a smooth 4-form on $\Gamma_{f,0}(r_0)$ it is 
straightforward that   $\int_{\Delta_{0'}(0,r)}dd^cT\to 0$.

\smallskip\noindent\sl
Step 2. \it  $\sup_{z'\in \Delta^{(n-1)}(0,r)}\int_{\Delta_{z'}^2(0,r)}dd^cT\to 0$ as $r\to 0$.

\smallskip\rm Othervice we would find a sequence $z_n'\to 0$ and $r_n\to 0$ such that 
$\int_{\Delta_{z_n'}(0,r_n)}dd^cT\le \eps_0<0$. Take any $r_0>\rho >0$ such that 
$\partial\Delta^2_{(0',\rho )}\cap 
{\cal I}(f)=\emptyset $ and remark that $\partial\Delta^2_{(z',\rho )}\cap 
{\cal I}(f)=\emptyset $ for  $z'$ close to zero. Then 

$$
\eps_0\ge \int_{\Delta_{z_n'}^2(0,r_n)}dd^cT\ge  \int_{\Delta_{z_n'}^2(0,\rho )}dd^cT=
$$

$$
=\int_{\partial \Delta_{z_n'}^2(0,\rho )}d^cT\to \int_{\partial \Delta_{0'}^2(0,\rho )}d^cT=
\int_{\Delta_{0'}^2(0,\rho )}dd^cT.
$$
But as $\rho \to 0$ the last integral tends to zero by step 1. This contradiction proves the lemma.

\smallskip\hfill{q.e.d.}

\smallskip\noindent\sl
2.4. \sl Prove in higher dimensions: pluriclosed case. 

\smallskip\rm
Fix a point $a\in A$ and suppose that there is a transversal sphere $\ss^3 = \{ x\in P: 
\Vert x-a\Vert =\eps \} $ on some two-plane through $a$ such that $f(\ss^3)$ is 
homologous to zero in $X$. We shall prove that in this case $f$ meromorphically extends 
to the neighborhood of $a$. Write $W=B^{n-1}
\times B^2$ for some neighborhood of this point $a\in A$ such that $(\bar B^{n+1}\times \partial
B^2)\cap A=\emptyset $ and for every  $z^{'}\in B^{n-1}$ one  has $f(\partial B^2_{z'})\sim 0$.
First we prove the following:

\smallskip\noindent\bf
Lemma 2.8. \it Suppose that the metric form $w$ on $X$ is pluriclosed and for
 all $z'\in B^{n-1}$ $f(\partial B^2_{z'})\sim 0$ in $X$. Then:

(i) $dd^c\tilde T=0$ in the sense of distributions.

(ii) There exists a $(1,0)$-current $\gamma $ in $W$, smooth in $W\setminus A$,
such that $\tilde T=i(\partial\bar\gamma -\bar\partial\gamma )$.
\smallskip\noindent\sl
Proof. \rm (i) Let $\tilde T_{\eps }$ be  smoothings of $\tilde T$ by
convolution. Then $\tilde T_{\eps }$ are plurinegative and $\tilde T_{\eps }
\to \tilde T$ in ${\cal D}_{n,n}(W)$. We have that
$$
\int_Wdd^c\tilde T_{\eps }\wedge (dd^c\Vert z'\Vert^2)^{n-1}=
\int_{\partial W}d^c\tilde T_{\eps }\wedge (dd^c\Vert z'\Vert^2)^{n-1} =
$$
$$
= \int_{\partial B^{n-1}\times B^2}d^c\tilde T_{\eps }\wedge (dd^c\Vert z'\Vert^2)
^{n-1} + \int_{B^{n-1}\times \partial B^2}d^c\tilde T_{\eps }\wedge (dd^c\Vert z'\Vert
^2)^{n-1}. \eqno(2.4.1)
$$
The first integral vanishes by  degree considerations. So
$$
\Vert dd^c\tilde T_{\eps }\wedge (dd^c\Vert z'\Vert^2)^{n-1}\Vert (W) =
-\int_{B^{n-1}}(dd^c\Vert z'\Vert^2)^{n-1}\int_{\partial B^2_{z'}}d^c\tilde T_{
\eps }. \eqno(2.4.2)
$$
Observe now that $\int_{\partial B^2_{z'}}d^c\tilde T_{\eps }\to
\int_{\partial B^2_{z'}}d^c\tilde T=\int_{f(\partial B^2_{z'})}d^cw=0$,
because $f(\partial B^2_{z'})\sim 0$ in $X$. So the right hand side of (2.4.2)
tends to zero as $\eps \searrow 0$. We obtain that
$$
\Vert dd^c\tilde T\wedge (dd^c\Vert z'\Vert^2)^{n-1}\Vert (W) =
\lim_{\eps \searrow 0}\Vert dd^c\tilde T_{\eps }\wedge (dd^c\Vert z'\Vert^2)
^{n-1}\Vert (W)=0.\eqno(2.4.3)
$$
\noindent
Taking sufficiently many such coordinate systems we see that $\Vert dd^c\tilde
T\Vert (W)=0$.
\smallskip\noindent
(ii) $\partial \tilde T$ is a $\bar\partial $-closed and $\partial $-closed
$(2,1)$-current. So, if $\phi \in {\cal D}_{n-1,n+1}(W)$ is $\partial $-closed
and such that $\bar\partial \phi = \partial\tilde T$ then $\phi $ is smooth on
$W\setminus A$ by elliptic regularity of $\bar\partial $. We have now
$d\tilde T=\partial \tilde T+\bar\partial \tilde T=\bar\partial \phi+\partial
\bar\phi $. Thus $d(\tilde T-\phi -\bar\phi )=0$. So $\tilde T-\phi -\bar\phi $
 is a $d$-closed current of degree two on $W$. Consider the following elliptic 
 system in
$W$:
$$
d\gamma =\tilde T-\phi -\bar\phi \hbox{ , }  d^*\gamma =0 \eqno(2.4.4)
$$
\noindent
Then (2.4.4) has a solution in $W$. Indeed, let $\gamma_1$ be any solution
of the first equation. Find a distribution $\eta $ on $W$ with $*d*d\eta
=\Delta \eta = *d*\gamma_1$ and put $\gamma_2 =\gamma_1-d\eta $.
$\gamma_2 $
is smooth on $W\setminus A$ because $\Delta \gamma_2 = d^*d\gamma_2 +
dd^*\gamma_2 = d^*(\tilde T-\phi - \bar\phi )$.
Write $\gamma_2 =i(\gamma^{1,0}-\bar\gamma^{1,0})$ -- the general form of a 
real 
$1$-form. We have
$i\partial \gamma^{1,0}=-\phi $ and $i\bar\partial
\bar\gamma ^{1,0}=+\bar\phi $, so
$$
\tilde T = d\gamma_2 +\phi + \bar\phi = d(i\gamma^{1,0}-i\bar\gamma^{1,0}) 
-i\partial \gamma^{1,0}+i\bar\partial \bar\gamma^{1,0} = i(-\partial \bar
\gamma ^{1,0} +
\bar\partial \gamma^{1,0}) \eqno(2.4.5)
$$
\noindent
with $\gamma^{1,0}$ having the required regularity. Now $\gamma =-\gamma^{1,0}$ 
satisfies (ii).
\smallskip
\hfill{q.e.d.}
\bigskip\noindent\bf
Lemma 2.9. \it If $\tilde T$ is pluriclosed, then the volumes $\Gamma_{f_z'}
\cap B^2_z\times X$ are uniformly bounded for $z\in B^{n-1}_r$ and $f$ 
extends meromorphically onto $W$.
\smallskip\noindent\sl
Proof. \rm Set $S:=T+dd^c\Vert z\Vert^2$, where $z=(z_n,z_{n+1})$. Note that 
$S$ is pluriclosed if $T$ is. Find $\gamma^{1,0}$ for $S$ as in Lemma 2.8, i.e.,
$\gamma^{1,0}$ is a $(0,1)$-current on $W$, smooth on $W\setminus A$, s.t. 
$S=i(\partial \bar\gamma^{1,0} - \bar\partial \gamma^{1,0})$. Smoothing by convolutions we 
still have
$\tilde S_{\eps }=i(\partial \bar \gamma^{1,0}_{\eps }-\bar\partial \gamma^{1,0}
_{\eps })$. Then for $z'\in B^{n-1}$ and $A_{z'}:=A\cap B^2_{z'}$ we have:

$$
{\sl vol}(\Gamma_{f_{z'}})=\int_{B^2_{z'}\setminus A_{z'}} S^2 =
\lim_{\eps \searrow 0}\int_{B^2_{z'}\setminus A_{z'}}\tilde S^2_{\eps }\le
\lim_{\eps \searrow 0}\int_{B^2_{z'}}\tilde S^2_{\eps }=
\lim_{\eps \searrow 0}\int_{B^2_{z'}}
i^2(\partial \bar\gamma^{1,0}_{\eps }-\bar\partial \gamma^{1,0}_{\eps })^2
\le
$$
$$
\le \lim_{\eps \searrow 0}\int_{B^2_{z'}}i^2d(\bar\gamma^{1,0}_{\eps }-
\gamma^{1,0}_{\eps })\wedge
d(\bar\gamma^{1,0}_{\eps }-\gamma^{1,0}_{\eps }) =
\lim_{\eps \searrow 0}\int_{\d B^2_{z'}}i^2(\bar\gamma^{1,0}_{\eps }-
 \gamma^{1,0}_{\eps })\wedge
d(\bar\gamma^{1,0}_{\eps }-\gamma^{1,0}_{\eps })=
$$
$$
= \int_{\partial B^2_{z'}}i^2(\bar\gamma^{1,0}-\gamma^{1,0})\wedge
d(\bar\gamma^{1,0}-\gamma^{1,0})\le const. \eqno(2.4.5)
$$
In the first inequality we  used the positivity of $T$. In the second -- the
fact that $-i^2\bar\partial \bar\gamma_{\eps }^{1,0}\wedge \partial \gamma_{\eps }
^{1,0}$ is positive and $\bar\partial \bar\gamma_{\eps }^{1,0}\wedge \bar
\partial \bar\gamma_{\eps }^{1,0}=0$. Finally $\gamma_{\eps }^{1,0}\to \gamma
^{1,0}$ on $\bar B^{n-1}\times \partial B^2$, since $\gamma ^{1,0}$ is
smooth there. This gives the required bound for  ${\sl vol}(\Gamma_{f_{z'}})=
\int_{B^2_{z'}\setminus A_{z'}}S^2$.

Lemma 1.3 (with $k=2$) gives us now the extension of $f$ onto $W\cong B^{n-1}\times B^2$.

\smallskip
\hfill{q.e.d.}

We end up with two remarks about the structure of the singularity set $A$ of 
our mapping in the presence of a pluriclosed metric form.

\smallskip
Consider two natural projections $\pi^1 :\Delta^{n+1}\to \Delta^{n-1}
\times \Delta_{z_n}$ and  $\pi^2 :\Delta^{n+1}\to \Delta^{n-1}
\times \Delta_{z_{n+1}}$. Remark that $\pi^j\mid_A$ are proper, $j=1,2$.
Set $A_j=\pi^j(A)$. We shall prove that each $A_j$ 
is pseudoconcave in 
$\Delta^n$ and admits a Sadullaev potential. We start with

\smallskip\noindent\bf
Lemma 2.10. {\it $A_j$ are complete pluripolar and moreover admit a Sadullaev 
potential.
}
\smallskip\noindent\bf
Proof. \rm Recall, see [Lv-Sl],  that a Sadullaev potential for a closed 
complete pluripolar set $A_j\subset \Delta^n$ is a plurisubharmonic 
function $\psi_j$ in $\Delta^n$ such that $\psi_j$
 is pluriharmonic on  $\Delta^n\setminus A_j$ and $A_j=\{ z\in \Delta^n: \psi_j(z)=
-\infty \} $.

For $z'=(z_1,...,z_{n})\in \Delta^n\setminus A_j$ define the area function 
$a_j(z')$ as in (2.2.3). The proof of Lemma 2.4 shows without essential changes 
that $dd^ca_j$ is a smooth $(1,1)$-form in $\Delta^n$. We claim that $a_j(z')\to 
+\infty $ as $z'\to A_j$. 

Suppose not, i.e. there is a sequence $\{ p_k\} \subset \Delta^n\setminus A_j$ 
such that $p_k\to p_0\in A_j$ and $\psi_j(p_k)\le C_0$ for some $C_0$ 
and all $k$. Corollary 1.7 in this case provides a complex curve $W$ in the 
neighborhood of $p_0$ such that $f$ meromorphically extends onto $W\times 
\Delta $. Using this it is straightforward to produce a transversal sphere 
$\ss^3\subset \Delta^{n+1}\setminus A$ such that $f(\ss^3)\sim 0$ in $X$. 
Contradiction with the Main Theorem. 

So $a_j(z')\to \infty $ as $z'\to A_j$. Find now some smooth function $h_j(z')$ 
such that $dd^ch_j=dd^ca_j$ and put $\psi_j=h_j-a_j$.

\smallskip\hfill{q.e.d.}

Recall that a closed subset $A_j\subset \Delta^n$ is pseudoconcave if $\Delta^n
\setminus A_j$ is pseudoconvex. 

\smallskip\noindent\bf
Corollary 2.11. {\it $A_j$ are pseudoconcave.}
\smallskip\noindent\bf
Proof. \rm It is sufficient to prove that for a holomorphic map $\phi :\Delta^2(r) 
\subset\Delta^n$ if $\phi :H^2(r)\to \Delta^n\setminus A_j$ then $\phi (\Delta^2)\subset 
\Delta^n\setminus A_j$. 

Let $\psi_j$ be a Sadullaev potential of $A_j$. Then $\psi_j\circ \phi $ 
is pluriharmonic on $H^2(r)$ and therefore extends pluriharmonically onto the 
$\Delta^2$. Therefore $\phi (\Delta^2)\cap A_j=\emptyset $.

\smallskip\hfill{q.e.d.}

\smallskip\noindent\bf
Remark. \rm In [Lv-Sl] an example of a closed subset $A\subset \Delta\times A(r,1)$
is constructed, which is pluripolar, pseudoconcave and admit a Sadullaev potential.
At the same time for every $z\in \Delta $ the set $A_z=\Delta_z\cap A$ is of Cantor-type
and $A$ admits no analytic structure.

\smallskip Some cases when $A$ possesses an analytic structure will be considered 
in \S \S 3.4, 3.5.

\smallskip\noindent\sl
2.5. The case of tamed structures.
\smallskip\rm
Let us now prove  Corollary 3 from the Introduction. Namely, we suppose that the
metric form $w$ on our space $X$ is the $(1,1)$-component of some closed
real two-form $w_0$, i.e., that there is a $(2,0)$-form $w^{2,0}$ such that
$w_0=w^{2,0}+w+\bar w^{2,0}$ and $dw_0=0$.

As we  remarked in the Introduction such $w$ is obviously $dd^c$-closed.
Thus the machinery of the proof of Main Theorem applies to this case. Therefore our
mapping $f$ can be extended meromorphically onto $\Delta^{n+1}\setminus A$,
where $A$ is either empty or is analytic of pure codimension two.

Suppose $A\not= \emptyset $. Take a point $a\in A$ with a neighborhood
$W\ni a$ biholomorphic to $B^{n-1}\times B^2$ and such that $\pi \mid_{\hat
A\cap W}:B^{n-1}\times B^2\to B^{n-1}$ is proper. Here $\hat A = A\cup I(f)
$ is a union of $A$ with the set of points of indeterminancy of $f$. Let us
prove that $dd^c\tilde T=0$ in $W$, where $T=f^*w$ on $\Delta^{n+1}\setminus
\hat A$.

From Lemma 2.15 we see that all we must prove is that
$\int_{\partial B^2_{z'}}d^c\tilde T_{\eps }=0$ for all $z'\in B^{n-1}$.
Indeed, let $T^0=f^*w_0$ and $T^{2,0}=f^*w^{2,0}$ on $\Delta^{n+1}
\setminus \hat A$. Then, since $dT^0=d^cT^0=0$, one has:

$$
\int_{\partial B^2_{z'}}d^c\tilde T_{\eps }=
\int_{\partial B^2_{z'}}d^c(\tilde T_{\eps }-T_{\eps }^0)=
\int_{\partial B^2_{z'}}d^c(-\tilde T_{\eps }^{2,0}-\bar T_{\eps }^{2,0}).
\eqno(2.5.1)
$$
Take a cut-off function $\eta $ with support in a neighborhood of
$ B^2_{z'}$. Then

$$
\int_{\partial B^2_{z'}}d^c\tilde T_{\eps }^{2,0}=
\int_{\partial B^2_{z'}}d^c(\eta \tilde T^{2,0})_{\eps }=
\int_{B^2_{z'}}dd^c(\eta \tilde T^{2,0})_{\eps }=0
\eqno(2.5.2)
$$
\noindent
by the reasons of bidegree.

So $\tilde T$ is pluriclosed on $W$ and we can extend $f$ onto the whole $W$ 
using Lemma 2.9.

\smallskip
\hfill{q.e.d.}

\bigskip\noindent\bf
3. Examples and open questions.

\smallskip\noindent\sl
3.1. General conjecture.
\smallskip\rm
Here we shall try to propose a general conjecture about extension 
properties of meromorphic mappings. In \S 1.5 we introduced the class ${\cal G}_k$ 
of reduced complex spaces possessing a strictly positive $dd^c$-closed $(k,k)$-form.

We conjecture that
meromorphic mappings into the spaces of class ${\cal G}_k$ are ''almost
Hartogs-extendable'' in bidimension $(n,k)$ for all $n\ge 1$. Namely, let 

$$
H^k_n(r):=\Delta^n_{1-r}\times \Delta^k\cup \Delta^n\times A^k(r,1)
$$

\noindent be the $k$-concave Hartogs figure.
\par\smallskip\noindent\sl
Cojecture. \it Every meromorphic map $f:H^k_n(r)\to X$, where $X\in
{\cal G}_k$ and is disk-convex in dimension $k$, extends to a meromorphic
map from $\Delta^{n+k}\setminus A$
to $X$,
where $A$ is a  closed $(n-1)$-polar subset of Hausdorff $(2n-1)$-dimensional
measure zero. Moreover, if $A\not=\emptyset $, then for every transversal 
sphere $\ss^{2k+1}$ in $\Delta^{n+k}\setminus A$ its image $f(\ss^{2k+1})$
is not homologous to zero in $X$. I.e. if $A\not= \emptyset $ then $X$ should 
contain a $(k+1)$-dimensional spherical shell. 
\par\smallskip\rm
In this paper we proved this conjecture in the case $k=1$.

\smallskip\noindent\sl
3.2. Examples to the Hartogs-type theorem.
\smallskip\rm 
 We start with an  example due to M. Kato, see [Ka-4].
\smallskip\noindent\bf
Example 3.1. \rm In $\cc\pp^3$ with homogeneous coordinates $[z_0:z_1:z_2:z_3]$ consider the
domain $D= \{ z\in \cc\pp^3: \vert z_0\vert^2+\vert z_1\vert^2>\vert z_2\vert
^2+\vert z_3\vert^2\} $. The natural action of $Sp(1,1)$ on $\cc\pp^3$ 
preserve $D$, i.e., $g(D)=D$ for all $g\in Sp(1,1)$. This action  is
transitive on $D$ and Kato proved, using  result of Vinberg, that there
exists a discrete subgroup $\Gamma\subset Sp(1,1)$, acting properly and
discontinuously on $D$, and such that $D/\Gamma =X^3$ is a compact complex
manifold, see [Ka-4] for details.

The projective plane $\cc\pp^2=\{ z_3=0\} $ intersects $D$ by the complement to
the closed unit ball $\bar B^4\subset \cc\pp^2$, namely, by
  $\cc\pp^2\setminus \bar B^4=\{ [z_0:z_1:z_2]\in \cc\pp^2:\vert z_2\vert^2<
  \vert z_0\vert^2+
\vert z_1\vert^2\} $. If $\pi :D\to X^3$ is the natural projection, then its
restriction $\pi\mid_{\cc\pp^2\cap D}:\cc\pp^2\setminus B^4\to X^3$
defines a holomorphic map from the complement to the closed unit ball to
$X^3$, which has a singularity at each point of $\partial B^4$!

However it is not difficult to see that  
$X$ doesn't contains neither two-dimensional spherical shells no 
three-dimensional ones.

This example shows also that one cannot hope for the better than it is 
conjectured in \S 3.1. Following two examples illustrate how the bounded 
cycle space geometry condition can be violated.

   \smallskip\noindent\bf
   Example 3.2. {\it There exists a holomorphic mapping 
   $f:\Delta \times \Delta_{{1\over 2}}\times A({1\over 2},1)\to X^3$ such that:

   (1) for any $s\in S=\{ (z_0,z_2)\in \Delta \times \Delta_{{1\over 2}}:
   \vert z_0\vert^2> \vert z_2\vert^2 \} $ the restriction
    $f_s=f\mid_{A_s(r,1)}$ extends holomorphically  onto $\Delta_s$;

   (2) for any $t>1$ there is a constant $C_t<\infty $ such that for all
   $s\in S_t=\{ (z_0,z_2)\in \Delta \times \Delta_{{1\over 2}}:
   \vert z_0\vert^2>t\cdot \vert z_2\vert^2 \} $ one has
   ${\sl area}(\Gamma_{f_s}) \le C_t$;

   (3) but for all $z\in \Delta^2\setminus \bar S =\{ (z_0,z_2)\in
   \Delta \times \Delta_{{1\over 2}}
   : \vert z_0\vert^2<
   \vert z_2\vert^2 \} $ the innner
  circle of the annulus $A^1_z(r,1):=\{ z_1\in \Delta_z: 1>\vert z_1\vert^2 >
  r^2 \} $ consists of essentially
  singular points of $f_z:A_z(r,1)\to X$, here $r^2=
  \vert z_2\vert^2  - \vert z_0\vert^2 $.
}

\smallskip\rm
Blow up $\cc\pp^4$ at the origin of its affine part. Denote by $\cc\pp^4_0$
the resulting manifold. There exists the  natural holomorphic projection
$p:\cc\pp^4_0\to \cc\pp^3$, where $\cc\pp^3$ is considered as the exceptional 
divisor. $\Gamma $
being a group of $4\times 4$ matrices acts naturally on the affine part
$\cc^4$ of $\cc\pp^4$. This action obviously extends onto $\cc\pp^4$ and
lifts onto $\cc\pp^4_0$. Moreover the actions of $\Gamma $ on $\cc\pp^3$
and $\cc\pp^4_0$ are equivariant with respect to the projection $p$.
Put $\hat D:=p^{-1}(D)$ and $\hat X:=\hat D/\Gamma $. Note that $p$ descends 
to a holomorphic map (in fact, $\cc\pp^1$-fibration) $p:\hat X\to X^3$. 
 If we take now
$\hat f:\{ z\in \cc^4: \vert z_0\vert^2+\vert z_1\vert^2>\vert z_2\vert^2,
z_3=0\} \to \hat X$ to be the restriction of quotient map, we get a mapping 
into $\hat X$ with properties (1),(2),(3). Taking composition $f:=p\circ \hat f:
\Delta\times \Delta_{{1\over 2}}\times A({1\over 2},1) \to X^3$ we get an 
example  of the 
same type with a $3$-dimensional image manifold.

\smallskip\noindent\bf 
Remark. \rm The rational cycle space in this example has a noncompact irreducible component 
of dimension four.

\smallskip For the detailed construction of the following example we refer 
the reader to [Iv-4], where another property of this example was studied.
Here we only list the properties of this example related to   
the cycle space geometry. Namely, denote 
by $z=(z_1,z_2)$ the coordinates in $\cc^2\times \{ 0\} \subset \cc^3$:
 
\smallskip\noindent\bf
Example 3.3. {\it There exists a compact complex manifold $X^3$ of dimension 
three and a meromorphic map $f:\Delta^3\setminus \{ 0\} \to X^3$ such that:

(1) for every cone $K_n:= \{ z=(z_1,z_2)\in \Delta^2: \vert z_2\vert >\vert z_1
\vert^n \}$ there is a constant $C_n$ such 

that ${\sl area}(\Gamma_{f_z})\le C_n$;

\smallskip\noindent but 

(2)  ${\sl area}(\Gamma_{f_z})\to \infty $ where $z=(z_1,0)$ and $z_1\to 0$;

(3) $f_0$ extends from $\Delta_0\setminus \{ 0\} $ onto $\Delta_0$;

(4) for every $t\in \cc\pp^1$ $\lim_{z\to 0}\Gamma_{f_z}$ (where 
$z=(z_1,tz_1^n)$) is equal to $\Gamma_{f_0}$ plus a rational 

cycle $Z_{n,t}$ which consists of $n$ components (counted with multiplicities);

(5) $\{ Z_{n,t}:t\in \cc\pp^1\} $ form an irreducible component $A_n$ of 
${\cal R}(X^3)$ and $\{ A_n\}_{n=1}^{\infty }$ is a 

connected chain of 
irreducible components of ${\cal R}(X^3)$.
  
}

\smallskip\noindent\bf
Remark. \rm We see that the space ${\cal R}(X)$ can contain an infinite connected 
chain of compact irreducible components, thou having non-bounded rational cycle geometry.
\smallskip\rm
We propose to the interested reader to check the details.  Of course, $f$ is the 
same map  that was constructed and studied in details in [Iv-4].

\medskip\rm
A statement similar to {\sl Lemma 1.3} (in fact a bit weaker one) was 
implicitly used in [Sb] Th\'eor\`eme 5.2, and was claimed to follow from Fubini's
theorem. The claim was that if one can bound the volumes of a meromorphic 
graph along sufficiently many two-dimensional directions then one can bound 
the total volume. If this were true, this would give a new proof of 
Siu's theorem of the removability of codimension two singularities for 
meromorphic mappings into  compact K\"ahler manifolds. It would
also replace our  Lemma 1.3 in the proof of the Main Theorem.

However, statements of such type cannot be derived from Fubini's theorem,
because the measure on the graph is not the product measure of the measures 
on the slices. Moreover, let
us see that an analogous statement is not valid in the real case (and thus
there is no "formula" to prove such things).

Of course, it should be made clear that this small gap doesn't touch any of 
the main results of [Sb].

\smallskip\noindent\bf
Example 3.4. {\it There exists a sequence of smooth mappings $f_n$ from the 
square  $\Pi = [0,1]\times [0,1]\subset \rr^2$ to $\rr^6$ such that:

(a) the lengs of the curves $f_n(x,\cdot ):[0,1]\to \rr^6$ and 
$f_n(\cdot ,y):[0,1]\to \rr^6$ are uniformly bounded for all $x,y\in [0,1]$,
but 

(b) the areas of $f_n(\Pi )$ turn to infinity.
}

\smallskip\rm

First of all one easily constructs a sequence of smooth strictly positive
functions $\phi_n$ on the square $\Pi = [0,1]\times [0,1]$, which:

a) are equal to $1$ in some fixed neighborhood of $\partial \Pi$;

b) for all $x,y\in [0,1]$

$$
\int_0^1\phi_n(x,t)dt\le 2 \hbox{ and } \int_0^1\phi_n(t,y)dt\le 2;
$$

c) $ \int_0^1\int_0^1\phi_n^2(x,y)dt\ge n$.

\smallskip
Consider  Riemannian metrics
$ds_n^2=\phi_n^2dx\otimes dx + \phi_n^2dy\otimes dy$ on the square. The length
of any segment parallel to the axis in this metric is
$$
\int_0^1\sqrt{\phi_n^2(x,t)+\phi_n^2(x,t)}dt\le 2\sqrt{2},
$$
while the area is $ \int_0^1\int_0^1\sqrt{\phi_n^2(x,y)\cdot \phi_n^2(x,y)}dt
\ge n$.

Now one can isometrically imbed $(\Pi ,ds_n^2)$ into $\rr^6$, see [G].
This gives the required example.

\smallskip
The author should say at this point that he doesn't know such an example 
with meromorphic $f_n$-s (and thus may be the implicit statement from [Sb] 
is nevertheless correct), and moreover he thinks that in  complex 
analytic setting such type of "area-volume" estimates could be true. 
The proof might follow from cycle space techniques like those used in the
present paper.

\smallskip\noindent\sl
3.3. Meromorphic correspondences.
\smallskip\rm
Let $D$ be a domain in complex space $\Omega$ and $x_0\in \partial D$ be a
boundary point. $D$ is said to be $q$-concave at $x_0$ if there is a
neighborhood $U\supset x_0$ and a smooth function $\rho :U\to \rr $ such that

1) $D\cap U  = \{ x\in U:\rho (x)<0\} $;

2) the Levi form of $\rho $ at $x_0$ has at least $n-q+1$ negative eigenvalues.
\smallskip\noindent
Here $n={\sl dim} \Omega $. By the Projection Lemma of Siu, see [Si-T], if
$x_0$ is a $q$-concave boundary point of $D$, $q\le n-1$, one can find 
neighborhoods $U\ni x_0$ and $V\ni 0\in \cc^n$ and a proper holomorphic map
$\pi :(U,x_0)\to (V,0)$ such that $\pi (D\cap V)$ will contain a Hartogs
figure $H$, whose associated polydisk $P$ contains the origin. Let $d$ be
the branching number of $\pi $.

Now suppose that a meromorphic map $f:D\to X$ is given, where $X$ is another
complex space. $f\circ \pi^{-1}$ defines a $d$-valued meromorphic
correspondence between $H$ and $X$.

\smallskip\noindent\bf
Definition 3.1. \it A $d$-valued meromorphic correspondence between complex
spaces $H$ and $X$ is an irreducible analytic subset $Z\subset H\times X$ such
that the restriction $p_1\mid_Z$ of the natural projection onto the first factor
on $Z$ is proper, surjective and generically $d$ - to - one.

\smallskip\rm
Thus the extension of $f$ onto the neighborhood of $x_0$ is equivalent to the
 extension of $Z$ from $H$ to $P$. It is clear that if $f$ was also a correspondence
it will produce no additional complications. Thus we should discuss how far the
problem of extending of correspondences go from the extension of mappings.

Let $Z$ be a $d$-valued meromorphic correspondence between the Hartogs figure
$H$ and $X$. $Z$ defines in a natural way a mapping $f_Z:H\to {\sl Sym}^d(X)$
- the symmetric power of $X$ of degree $d$. Clearly the extension of $Z$ onto $P$
is equivalent to the extension of $f_Z$ onto $P$. If $X$ was, for example, a
K\"ahler manifold, then  ${\sl Sym}^d(X)$ is a K\"ahler space by [V]. So,
meromorphic correspondences with values in K\"ahler manifolds are extendable
through pseudo-concave  boundary points.

For the manifolds from class ${\cal G}_1$ this is no longer the case,  
even if
they do not contain spherical shells.

\smallskip\noindent\bf
Example 3.5. \it There exists a compact complex (elliptic) surface $X$ such that:

(a) every meromorphic map $f:H^2(r)\to X$ extends meromorphically onto
$\Delta^2$, but

(b) there exists a two-valued meromorphic correspondence $Z$ between
$\cc^2_*$ and $X$ 
that cannot be extended  to the origin.
\smallskip\rm
Consider the standard Hopf surface $H=\cc^2\setminus \{ 0\} /(z\sim 2\cdot z)$.
Denote by $\pi :H\to \cc\pp^1$  the standard projection.
Let $\phi : C\to \cc\pp^1$ be a nonconstant meromorphic function on the
Riemann surface $C$ of positive genus. $\phi $ is a $d$-sheeted
ramified covering of $\cc\pp^1$ by $C$. If we take $C$ to be a torus we can
have such $\phi $ with $d=2$. Following Kodaira we shall construct
an elliptic surface over $C$ in the following way. Put
$$
X_1 = \{ (z,y)\in C\times H : \phi (z)=\pi (y)\}.\eqno(3.3.1)
$$
\noindent
Elliptic structure on $X_1$ is given by the restriction onto $X_1$ of the
natural projection $p_1:C\times H\to C$. Note that the restriction 
of the natural projection $p_2:C\times H\to H$ onto $X_1$ gives us a $d$-sheeted
covering $p_2\mid_{X_1}$ of $H$ by $X_1$ preserving the elliptic
structure. Let $n:X\to X_1$ be the normalization of $X_1$. Then $X$ is a
smooth elliptic surface over $C$ with elliptic fibration $p:=p_1\mid_{X_1}
\circ n:X\to C$, and $F:=p_2\mid_{X_1}\circ n :X\to H$ is a  $d$ -
sheeted covering.

$Z:=F^{-1}\circ \pi : \cc^2_*\to X$ is the $d$-valued meromorphic correspondence
between  $\cc^2_*$ and $X$, which cannot be extended to the origin, because
the projection $\pi :\cc^2_*\to H$ cannot be extended meromorphically to
zero.

On the other hand in [Iv-1] it was proved that meromorphic mappings from
$H^2(r)$ to an elliptic surface over a Riemann surface of positive genus
are extendable onto $\Delta^2$.
\smallskip
One can interpret this example in the way that ${\sl Sym}^2(X)$ may have
a spherical shell  even when $X$ has no.

\smallskip\rm
In view of discussion above we can restate our results for
meromorphic correspondences.
\smallskip\noindent\bf
Definition 3.2. \it By a branched spherical shell of degree $d$  in a
complex space $X$ we shall understand the image $\Sigma $ of $\ss^3\subset \cc^2
$ under a $d$-valued meromorphic correspondence between some neighborhood
of $\ss^3$ and $X$ such that $\Sigma \not\sim 0$ in $X$.

\smallskip\noindent\bf
Corollary 3.1. \it Let $Z$ be a meromorphic correspondence from the domain
$D$ in a complex space $\Omega $ into a  disk-convex
complex space $X\in {\cal G}_1$ and let $x_0$ be a concave boundary point of
$D$. Then $Z$  extends onto some neighborhood of $x_0$ in $\Omega $ minus
(possibly empty) complex variety $A$ of pure codimension two. If $X$ doesn't
contain  branched
spherical shells then $A=\emptyset $.
\smallskip\noindent\bf
Remark. \rm We would like to point  out here that a branched shell could be a
 much violate object than a non branched one. To see this 
consider a smooth
complex curve $C$ in a neighborhood of $\bar \bb^2(2)\setminus \bb^2(1/2)$
that doesn't extend to $\bb^2(1)$. Let $\pi :W\to \bb^2(2)\setminus \bar\bb^2
(1/2)$ be a covering branched along $C$. The images of $\hat S:=\pi^{-1}(\ss^3)
$ in  complex spaces could be  branched shells that would not bound
an abstract Stein domains.

\smallskip\noindent\sl 
3.4.  Mappings into compact complex surfaces.

\smallskip\rm Here we shall prove the Corollary 1 from the Introduction. Let $X$ denote 
now a compact complex surface. First of all we remark that  if $X$ is K\"ahler then by 
[Iv-3] every meromorphic map from $H^{n+1}_U(r)$ to $X$ extends onto $\Delta^{n+1}$. In 
particular the singularity set $A$ is empty. From Enriques-Kodaira classification we 
are left with two classes: elliptic surfaces and surfaces of class VII withought 
meromorphic functions. Since every compact complex surface carries a pluriclosed metric 
form the main theorem applies.  

\smallskip\noindent\sl Elliptic surfaces. \rm Let our map be extended onto $\Delta^{n+1}
\setminus A$ as in the  Main theorem. Take a point $a\in A$ and choose coordinates 
in a neighborhood $V$ of $a$ such that $V\cong \Delta^{n-1}\times \Delta^2$ and $f$ is 
holomorphic in the neighborhood of $\Delta^{n-1}\times \partial\Delta^2$. Now for 
every $a\in \Delta^{n-1}$ the $2$-disk $\Delta^2_a:=\{ a\} \times \Delta^2$ 
is "transversal" to $A$, i.e. $A\cap \Delta^2_a$ is compact of dimension  zero.
Consider the restriction $f_a$ 
of $f$ onto $B^2_a$. $f_a$ is defined and meromorphic outside of a compact set. Lemma 5 
from [Iv-1] shows that $f_a$ extends to the complement of finitely many points $\{ b_k(a)\} $. 
Moreover, the set $A_a:= \{ b_k(a)\} $ is not empty only in the case when $X$ is an elliptic 
surface over $\cc\pp^1$ with projection $\pi :X\to \cc\pp^1$ and $\{ b_k(a)\} $ is exactly 
the set of points of indeterminancy $I(\pi\circ f_a)$ of $\pi \circ f_a$. Since $I(\pi\circ
f_a)\subset I(\pi\circ f)$ and the 
latter is analytic, the Corollary 1 for elliptic surfaces follows.

\smallskip\noindent\sl Surfaces of class VII. \rm Toruses and K3 surfaces are K\"ahler, and therefore
for them $A$ is empty as above. We shall crusially use that for all others surfaces from the 
class VII $b_1=1$. First we prove  
that any meromoprhic map $f:\bb^2(1)\setminus \bb^2(r)\to X$ from the spherical shell into $X$ 
extends onto $\bb^2\setminus \{ $ finite set $\} $ and the number of points in this set is 
bounded by some number $N$. We shall show that $N\le \vert \int_{f(\partial B^2)}d^c\omega 
\vert / \vert \int_{M}d^c\omega \vert $, where 
 $M$ is the generator of the torsion free part of $H_3(X,\zz )$. 
 
 By the main theorem it extends 
 onto $\bb^2$ 
minus a compact $A$ of Hausdorff dimension zero. For any sphere $S$ around $b\in A$ s.t. 
$S\cap A=\emptyset $ its 
image $f(S)$ is not homologous to zero. Since the first Betti number of $X$ is one by the 
Poincar\'e duality the third is also one. Let $\eps_1=\int_ Md^c\omega $, where $\omega $ is a pluriclosed metric form on $X$.
Then for 
every $S$ there is a nonzero integer $n$ such that $f(S)$ is homologous to $nM$ modulo 
torsion.  
 The image $f(S)$ is not contained in the torsion 
part of $H_3(X,\zz )$ because $\int_{f(S)}d^c\omega \not= 0$. Therefore 
$\int_{f(S)}d^c\omega = n\eps_1$.

\smallskip\noindent\sl Case 1. $\eps_1=0$. 

\smallskip\rm This implies that all $\int_{f(S)}d^c\omega $ 
are equal to zero and by Lemma 2.5 $f$ extends meromorphically onto every such $b$. This means that 
$A$ is empty in this case.

\smallskip\noindent\sl Case 2. $\eps_1\not= 0$.   

\smallskip\rm Since $\int_{f(S)}d^c\omega =
n\eps_1$ for nonzero integer $n$ this implies that $\vert \int_{f(S)}d^c\omega \vert $ are 
separate from zero. From the other hand $\Sigma_{b\in A}\int_{f(S(b))}d^c\omega  = \int_{\partial \bb^2}
d^c\omega $ 
and therefore is finite. Therefore the set  $A$ is finite and $\vert A\vert \le 
\vert \int_{f(\partial B^2)}d^c\omega \vert / \vert \int_{M}d^c\omega \vert $.

\smallskip In the case of a higher dimension take a point $a\in A$ and choose a 
neighborhood $V\cong \Delta^{n-1}\times \Delta^2$ of $a$ such that $A\cap (\Delta^{n-1}
\times \partial \Delta^2)=\emptyset $. $A\cap V$ is a graph of an 
$N$-
valued continuous mapping of $\Delta^{n-1}$ to $\Delta^2$. $N$- valued means that
$\vert A_{z'}\vert \le N$ for every $z'\in \bar \Delta^{n-1}$ and there exists
$z_0^{'}$ such that $\vert A^0_{z'}\vert =N$. A multivalued mapping is 
continuous  if the set $A$ is
closed, which is obviously our casae. Note that the trivial extension $\tilde T$ of 
$T=f^*\omega $ is an $L^1$-current on 
$V$ with $dd^c\tilde T\le 0$  supported on $A$.

\smallskip\noindent\bf
Lemma 3.2. \it Let $A$ be the graph of an $N$- valued continuous mapping
of $\bar\Delta^k$ to $\Delta^l$ and let $R$ be a closed positive current in
$\Delta^{k+l}$ of bidimension $(k,k)$ supported on $A$. Then $A$ is a pure
 $k$-dimensional analytic variety in $\Delta^{k+l}$.
\smallskip\noindent\sl
Proof. \rm Write $R=R_{K,\bar J}({i\over 2})^k{\partial \over \partial z^K}
\wedge {\partial \over \partial \bar z^J}$, where $K$ and $J$ are multiindices
of length $k$. Consider the measures $R_{K,\bar J}$, denote by 
$\mu_{K,\bar J}=\pi_*(R_{K,\bar J})$ their direct images and desintegrate this 
measures
 with respect to the natural projection $\pi :\Delta^k\times \Delta^l\to
\Delta^k$. Desintegration means that one has probability measures
$\nu_{K,\bar J,z'}$ on $\Delta^l_{z'}:=\{ z'\} \times \Delta^l$ with the
property that for every continuous function $h$ in $\Delta^{k+l}$
$$
<R_{K,\bar J},h> = \int_{\Delta^k}(\int_{\Delta^l_{z'}}\bar h\mid_{\Delta^
l_{z'}}d\nu_{K,\bar J,z'})d\mu_{K,\bar J}, \eqno(3.4.1)
$$
\noindent
see [D-M].

Let $\Omega $ be the maximal open subset of $\Delta^k$ such that the
multivalued map $s$, which is given by its graph $A$, takes on exactly $N$
different values (and $N$ is maximal). First we shall prove that $A\cap
(\Omega \times \Delta^l)$ is analytic.

Let further $\Omega_1$ be some simply connected open subset of $\Omega $.
 Then $s\mid_{\Omega_1}$ decomposes to $N$ well-defined single valued maps $s^1,
...,s^N$. So it is enough to consider the case when $s$ is single valued. Put
$s(z')=(s_1(z'),...,s_l(z'))$. Note that in this case $\nu_{K,\bar J,z'}=
\delta_{\{ z^{''}-s(z')\} }$.  Therefore for the coefficients $R_{K,\bar J}$ 
of our current $R$
and for $\phi \in C^{\infty }(\Omega_1\times \Delta^{l})$ such that $\pi ({\sl supp} \phi
)\subset \subset \Omega_1$ we can write that
$$
<R_{K,\bar J},\phi > = \int_{\Omega_1 }\bar \phi (z',s(z'))d\mu_{K,\bar J}(z').
\eqno(3.4.2)
$$
If we choose $\phi $ not depending on $z^{''}:=(z_{k+1},...,z_{k+l})$ then
(3.4.2) gives
$$
<R_{K,\bar J},\phi > = \int_{\Omega_1}\bar\phi (z')d\mu_{K,\bar J}(z').
\eqno(3.4.3)
$$
From the closedness of $R$ we obtain that
$$
0=<R,d[({i\over 2})^k\bar\phi (z')
dz_1\wedge ...dz_{p-1}\wedge dz_{p+1}\wedge ...\wedge dz_k\wedge d\bar z_J]>
=<R_{1...k,\bar J},{\partial \bar\phi \over \partial z_p}> =
$$
$$
= \int_{\Omega_1}{\partial \phi \over \partial \bar z_p}d\mu_{1...k,\bar J}.
$$

So $\mu_{1...k,\bar J}(z') = c_{1...k,\bar J}(z')\cdot ({i\over 2})^kdz'\wedge d\bar
z'$, where $c_{1...k,\bar J}$ are holomorphic for all $J$. In
particular $c_{1...k,\bar 1...\bar k}$ is constant. Now take the $(k-1,k)$-
forms $\psi_{q\bar p} = \bar\phi (z')\cdot \bar z_p\cdot ({1\over 2})^kdz_1\wedge
...\wedge dz_{q-1}\wedge dz_{q+1}\wedge ...\wedge dz_k\wedge d\bar z_1\wedge ...\wedge 
d\bar z_k$. We
have
$$
0 = <R,d\psi_{q,\bar p}> = <R,{\partial \bar\phi \over \partial z_q}\cdot \bar z_p
({i\over 2})^kdz'\wedge d\bar z'> = <R_{1...k,\bar 1...\bar k},{\partial \bar\phi \over
\partial z_q}\bar z_{\bar p}> =
$$
$$
 = c_{1...k,\bar 1...\bar k} \int_{\Omega_1}{\partial \phi \over \partial \bar z_q}(z')
\cdot s_p(z')({i\over 2})^kdz'\wedge d\bar z',
$$
\noindent
i.e., $s_p$ are holomorphic.

Thus we have proved that $s$ is an $N$-valued analytic map of $\Omega $ into
$\Delta^l$. Considering appropriate discriminants and using Rado's Theorem,
we obtain analyticity of $s$ on the whole $\Delta^k$.

\smallskip
\hfill{q.e.d.}

\smallskip\noindent\sl
3.5. More remarks about the structure of the singularity set.\rm

 One can slightly generalise considerations above and prove in several other cases 
 that the singularity set $A$ of 
 a meromorphic map into a complex manifold with pluriclosed metrik form is analytic of 
 pure codimension two. We shall prove the following
 
 \smallskip\noindent\bf
 Theorem 3.3. {\it Suppose that a disk-convex complex spave $X$ admits a pluriclosed 
 metrik form $\omega $ such that $d^c\omega \in H^3(X,\zz )$. Then every meromorphic map 
 $f:H^{n+1}_U(r,1)\to X$ extends onto $\Delta^{n+1}\setminus A$, where $A$ is an analytic 
 subset of pure codimension two.
 }
 
 \smallskip\rm We start with the case $n=1$.

\smallskip\noindent\bf 
Lemma 3.4. {\it For $X$ as in Theorem 3.3  $f:H^2_U(r,1)\to X$ can be 
meromorphically extended onto $\Delta^2\setminus A$, where $A$ is discrete.
}
\smallskip\noindent\bf
Proof. \rm By the Main Theorem we know that $f$ extends onto $\Delta^2\setminus A$, 
where $A$ is pluripolar of Hausdorff dimension zero. We suppose that $A$ is a minimal subset 
of $\Delta^2$ such that $f$ extends onto $\Delta^2\setminus A$. Take a relatively compact 
open subset $P\subset \Delta^2$ such that $\partial P\cap A=\emptyset $  and choose a 
finite subcomplex $K$ of CW-complex $X$ to contain the 
$f(\bar P\setminus A)$, which is a compact subset of $X$ due to the disk convexity 
of the latter. All we need to prove is that $P\cap A$ is finite. 
Let $\theta_1,...,\theta_N$ be the generators
of the free part of $H^3(K,\zz )$ and $\psi_1,...,\psi_L$ be the generators of the 
free part of $H_3(K,\zz )$. 
Take integer numbers $r_1,...,r_N$ such that 

$$
d^c\omega  = r_1\theta_1 +...+r_N\theta_N. \eqno(3.5.1)
$$

Take a ball $B\subset\subset \Delta^2$ with $\partial B\cap A=\emptyset $. 
Then there are integers $z_1,...,z_L$ such that 
$$
f(\partial B) = z_1\psi_1+...+z_L\psi_L \eqno(3.5.2)
$$
\noindent 
in $H_3(K,\zz )$ modulo torsion. For the measure $\mu $ defined from $dd^c\tilde T=\mu \cdot 
\omega_e^2$ we have that 
$$
\mu (B\cap A)= \int_Bdd^c\tilde T = \int_{\partial B}d^cT.\eqno(3.5.3)
$$ 
Using that, we can write
$$
\mu (B\cap A) = \int_{f(\partial B)}d^c\omega = \sum_{k=1}^N\sum_{i=1}^Lz_ir_k 
\int_{\psi_i}\theta_k .\eqno(3.5.4)
$$
\noindent
Put $c^{ik}=\int_{\psi_i}\theta_k\in \zz $. Now if we put $\tilde z^k = 
\sum_{i=1}^Lz_ic^{ik} \in \zz $ then 
$$
\mu (B\cap A) = \sum_{k=1}^N\tilde z^kr_k. \eqno(3.5.5)
$$

\noindent
The right hand side of (3.5.5) is a negative integer and therefore is separated from zero, i.e.
there exists an $\eps_0<0$ such that $\mu (B\cap A)<\eps_0 $ if $B\cap A\not =\emptyset $.
Therefore, if $A$ is not discrete there is a sequence $\{ s_n\} 
\subset A$ converging to $s_0\in A$, and therefore if we take nonintersecting $B_{\eps_n}(s_n)$
we obtain $\mu (S)\le \Sigma_n\mu (B_{\eps_n}(s_n))=-\infty $. Contradiction.

\smallskip
\hfill{q.e.d.}

Further increasing of the dimension is now an obvious 
repetion of the case of surfaces.

\smallskip\hfill{q.e.d.}

\smallskip\noindent
{\sl 3.6.  More examples of a singularity set.}
\smallskip

Arguments from the previous sections will not work if $d^c\omega $ is not in 
$H^3(X,\zz)$. For example if there is two 3-cycles (holomorphic images of $\ss^3\subset 
\cc^2$ ) $M_1, M_2$ such that $\int_{M_1}d^c\omega =1$ and say $\int_{M_1}d^c\omega $ 
is irrational, then there could be a sequence of points $\{ a_k\} \subset A\subset B^2$ 
such that a meromorphic mapping $f:B^2\setminus A\to X$ has the property that for 
sufficiently small spheres $S_k$ around $a_k$ $f(S_k)=n_kM_1+m_kM_2$ with $\{ m_k\} $ 
unbounded. Then one could have that $\int_{f(S_k)}\omega \to 0$ and arguments from 
prevoius two sections will fail. $A$ can be now the the set of accumulation points 
of the sequence  $\{ a_k\}$ and can be a Cantor-type set. In this section we give 
such examples, both with $A=\{ a_k\}$ and $A$ uncountable. In the Example 3.6 below 
a map will behave exactly as described above. Unfortunatly  it is 
not clear if that manifold $X$ can be endoved with a pluriclosed (or at least plurinegative) 
metric form.

\smallskip\noindent\bf
Example 3.6. {\it There exists a compact complex threefold $X$ with $H^3(X,\zz )$ 
generated by two cycles $M_1, M_2$, which are holomorphic images of $\ss^3\subset \cc^2$
(i.e. spherical shells) and a holomoprhic mapping $f:B^2\setminus A:\to X$ such that: 

\smallskip
a) $A$ is an uncountable Cantor-type complete pluripolar compact subset of $B^2$;

b) all points of $A$ are singular for $f$;

c) there is a dense in $A$ sequence $\{ a_k\} \subset A$ and spheres $S_k$ around $a_k$ 
such that $f(S_k)=n_kM_1+m_kM_2$ with $\{ (n_k,m_k)\} $ unbounded.

}

\smallskip\rm 

Take $\cc\pp^3$ with homogeneous coordinates $[z_0:z_1:z_2:z_3]$ and 
take two lines $l_1=\{ z_0=z_1=0\} $ and $h_1= \{ z_2=z_3=0\} $. Consider 
the following function $\phi_1 (z)={\vert z_0\vert^2+\vert z_1\vert^2\over 
\vert z_0\vert^2+\vert z_1\vert^2 +\vert z_2\vert^2+\vert z_3\vert^2 }$.
Remark that $l_1=\{ z:\phi_1 (z)=0\}$ and $h_1=\{ z:\phi_1 (z)=1\} $. For 
a complex number $\alpha $ such that $\vert \alpha \vert <1$ consider 
a domain $D_{\alpha }^1\subset \cc\pp^3$ defined by $D_{\alpha }^1 =
\{ z\in \cc\pp^3:\vert \alpha \vert <\phi_1 (z)<\vert {1\over \alpha } 
\vert \} $. The group of 
automorphisms of $\cc\pp^3$ generated by $g_1(z)= [{1\over \alpha^2}z_0:
 {1\over \alpha^2}z_1:z_2:z_3]$ has $D_{\alpha}^1$ as its fundamental 
 domain.
 
 Take new coordinates 
 
 $$
 w_0 = {1\over \sqrt 2}(z_0+z_2),  w_2 = {1\over \sqrt 2}(z_1+z_3),
 $$
 
 $$ 
 w_1 = {1\over \sqrt 2}(z_0-z_2),  w_3 = {1\over \sqrt 2 }(z_1-z_3),
 $$

 \noindent and repeat the above considerations: define lines $l_2=
 \{ w_0=w_1=0\} $ and 
 $h_2= \{ w_2=w_3=0\} $, function $\phi_2 (z)={\vert w_0\vert^2+\vert w_1
 \vert^2\over 
\vert w_0\vert^2+\vert w_1\vert^2 +\vert w_2\vert^2+\vert w_3\vert^2 }$, take 
the domain $D_{\alpha }^2 =\{ w\in \cc\pp^3:\vert \alpha \vert <\phi_2 (w)<\vert {1\over 
\alpha }\vert \}$, and finally an automorphism $g_2(w)= [{1\over \alpha^2}w_0:
 {1\over \alpha^2}w_1:w_2:w_3]$.
 
 Now the domain $D_{\alpha }:=D^1_{\alpha }\cap D^2_{\alpha }$ will be 
 a fundamental domain for the group $G$ of biholomoprhic automorphisms of 
 $\cc\pp^3$ generated by $g_1, g_2$. Set $\Omega_{\alpha } \bigcup_{g\in G}
 g(D_{\alpha })$. Complement to $\Omega_{\alpha }$ in $\cc\pp^3$ is a Cantor 
 set of complex lines $A_{\alpha}$. If $\vert \alpha \vert $ was taken suficiently small
 the set $A_{\alpha }$ will be $2-polar$ and of Hausdorff $3$-dimensional measure zero.
 
 Set $X=\Omega_{\alpha }/G$ - this is our manifold. $H_3(X,\zz )$ is generated
 by $3$-cycles $M_1=\{ z_3=0, \vert z_0\vert^2+\vert z_1\vert^2=\vert z_2\vert^2\}
 $ and  $M_2=\{ w_3=0, \vert w_0\vert^2+\vert w_1\vert^2=\vert w_2\vert^2 \}
 $. Both are holomorhic images of the standart sphere from $\cc^2$. 
 
 Take the plane $P=\{ z_3=0\} $ and the ball $B:=\{ \phi_1\mid_P(z)<\vert \alpha 
 \vert \} $ on this plane. The the natural map $f:B\setminus A\to X$ has 
 all needed properties, where $A=B\cap A_{\alpha }$.
 
 It should be noted that the author doesn't knows wether or not this $X$
 admits a plurinegative metric form. 
 
 \smallskip\noindent\bf 
 Remark. \rm Manifolds of this type where first constructed by M. Kato in dimension 
 3, see [Ka-2], [Ka-3] and by M. Nori in higher dimensions, see [No].
 
 \smallskip\rm Let us give one more example with an interesting singular set.

\smallskip\noindent\bf
Example 3.7. {\it There exists a compact complex manifold $X$ of dimension 
four and a holomorphic map $f:\bb^2_*\setminus \{ ({1\over 2^{k-1}}\cdot {3\over 4} ,0)\}_{k=1}^{\infty }\to X$ 
such that each point $a_k=({1\over 2^{k-1}}\cdot {3\over 4} ,0)\} $ 
is an essential singularity of $f$. This means that $f$ doesn't extend to 
the neighborhood of such $a_k$ even meromorphically. Moreover, the 
manifold $X$ carries a plurinegative metrik form.
}
\smallskip\noindent\sl
Construction of $X$. \rm Take a Hopf surface $H=\cc^2\setminus \{ 0\} /(z\sim 
2z)$. By $\pi_1 :\cc^2\setminus \{ 0\} \to H$ denote the canonical prejection. 
Fix the  point $a_1=({3\over 4},0)$ in $\bb^2_*:=\bb^2\setminus \{ 0\} $ and 
its image $b:=\pi_1 (a_1)$ in $H$. Recall, see [Gr-Ha] p.   , that on a 
compact complex surface $H$ with $h^{0,2}=0$ there exists a holomorphic rank 
two vector bundle $p:E\to H$ having a holomorphic section $s:H\to E$ with 
$b$ being the only its zero point. The second Betti number of the Hopf 
surface is zero and therefore from $b^2=h^{2,0}+h^{1,1}+h^{0,2}$ we see 
that in our case the condition $h^{0,2}=0$ is satisfied. 

Denote by $Z$ the zero section of $E$.
Let $g:E\setminus Z\to E\setminus Z$ be the multiplication by $2$. The desired 
four-manifold $X$ is the quotient of $E\setminus Z$ by the group of 
biholomorphisms $G:= \{ g^n\}_{n\in \zz }$. Denote by $\pi_2:E\setminus Z
\to X$ the canonical projection.

\smallskip\noindent\sl 
Construction of the map $f$. \rm Denote by $\Gamma_s$ the graph of the 
section $s$ in the total space of the bundle $E$ and by $\Gamma_s':=
\Gamma_s\cap (E\setminus Z)$. Our map $f$ is defined to be the composition 
$f:=\pi_2\circ s\circ \pi_1:\bb^2_*\setminus 
\{ ({1\over 2^{k-1}}\cdot {3\over 4} ,0)\}_{k=1}^{\infty }\to X$.

Since our map satisfies $f({1\over 2}z)=f(z)$ the set $\{ ({1\over 2^{k-1}} \cdot {3\over 4},0)\}_{k=1}^{\infty }$
is precisely the set of all points where our $f$ is not defined. On 
$\bb^2_*\setminus \{ ({1\over 2^{k-1}}\cdot {3\over 4} ,0)\}_{k=1}^{\infty }$ the map $f$ is well defined 
and holomorphic.

To see that $f$ doesn't extend meromorphically to the neighborhood of any $a_k$ it 
is sufficient to show that for a small sphere $S_k$ around $a_k$ its image $f(S_k)$ 
is not homologous to zero in $X$. Take $S_k$ to be the boundary of the Euklidean 
ball $B_k$ with center at $a_k$ small enough to satisfy the following two conditions:

\smallskip
1) $E$ is trivialisible in the neighborhood $U_k$ of $\pi_1(B_k)$;

2) $B_k$ doesn't contain any  $a_j$ for $j\not= k$ and moreover $B_k$ is contained 
in the spherical region $\{ z\in \cc^2: {1\over 2^k}< \Vert z\Vert <{1\over 2^{k-1}}\} $.

\smallskip Image $\pi_2(U_k\times \cc^2_*)$ is a product $U_k\times H$. Now it is clear 
that $f(S_k)$ is homologous to the generator of $H_3(U_k\times H,\zz )=H_3(H,\zz )=\zz $.

\smallskip\noindent\bf
Remark. \rm If any meromorphic mapping $f:H^2(r)\to X$ (with $X$ admitting a pluriclosed 
metric form) extends onto $\Delta^2$ minus a countable set (as in the example abowe),
then using the theorem of Nishino, see [Ni], one can prove that every meromorphic 
map  $f:H^n(r)\to X$ extends onto 
$\Delta^n\setminus A$, where $A$ is at most countable union of locally closed analytic 
sets of pure codimension two.

\magnification=\magstep1
\spaceskip=4pt plus3.5pt minus 1.5pt
\spaceskip=5pt plus4pt minus 2pt
\font\csc=cmcsc10
\font\tenmsb=msbm10
\def\rr{\hbox{\tenmsb R}}
\def\cc{\hbox{\tenmsb C}}
\newdimen\length
\newdimen\lleftskip
\lleftskip=2.5\parindent
\length=\hsize \advance\length-\lleftskip
\def\entry#1#2#3#4\par{\parshape=2  0pt  \hsize%
\lleftskip \length%
\noindent\hbox to \lleftskip%
{\bf[#1]\hfill}{\csc{#2 }}{\sl{#3}}#4%
\medskip
}
\ifx \twelvebf\undefined \font\twelvebf=cmbx12\fi

\bigskip\bigskip
\bigskip\bigskip
\centerline{\twelvebf References.}
\bigskip




\entry{Ba-1}{Barlet D.:}{Espace analytique reduit des cycles analytiques
complexes compacts d'un espace analytique complexe de dimension finie .}
 Seminar Norguet IX, Lect. Notes Math., {\bf 482}, 1-157, (1975).

\entry{Ba-2}{Barlet D.:}{Majorisation du volume des fibres g\'en\'eriques et la 
forme g\'eom\`etrique du th\'eor\`eme d'applatissement.} Seminar P. Lelong - H. Skoda,
Lecture Notes in Math., {\bf822}, 1-17, (1980).





\entry{D-M}{Dellacherie C., Meyer P.-A.:}{Probabilit\'es et potentiel.} Publ.
\ Inst.\ Math.\ Univ.\ Strasbourg, {\bf XV}, Hermann (1975).





\entry{Fd}{Federer :}{Geometric measure theory.} Berlin, \ Springer (1969).

\entry{F-Sb}{Fornaess J.-E., Sibony N.:}{Oka's unequality for currents and
applications.} \ Math.\ Ann . {\bf301}, N4, 813-820 (1995).

\entry{Fj}{Fujiki A.:}{Closedness of the Douady Space of Compact K\"ahler
Spaces.} Publ.\ RIMS,\ Kyoto \ Univ. {\bf14}, 1-52, (1978).

\entry{Ga}{Gauduchon P.:}{Les metriques standard d'une surface a premier
nombre de Betti  pair.} Asterisque.\ Soc.\ Math.\ France.  {\bf126}, 129-135,
(1985).

\entry{Gl}{Golysin G.:}{Geometric function theory.} Nauka,\ Moscow (1966).

\entry{Gr}{Griffiths P.:}{Two theorems on extensions of holomorphic mappings
.} Invent.\ math. {\bf14}, 27-62, (1971).

\entry{G}{Gromov M.:}{Partial Differential Relations.} Springer, \ Ergebnisse
der \ Mathematik und ihrer \ Grenzgebiete, 3 Folge, Band 9, (1986).

\entry{Ha}{Hartogs F.:}{Zur Theorie der analytischen
 Funktionen me\-hrerer un\-abh\"angiger Ver\"andenlichen insbesondere
\"uber die Darstellung derselben durch Reihen, wel\- che nach Potenzen
 einer Ver\"anderlichen fortschreiten.} \  Math.\ Ann. {\bf62}  1-88,
(1906).
  

\entry{H-S}{Harvey R., Shiffman B.:}{A characterization of holomorphic chains
.} Ann.\ Math. {\bf99}, 553-587, (1974).




\entry{Hm}{H\"ormander L.:}{Notions of convexity.} Birkh\"auser (1994).

\entry{Iv-1}{Ivashkovich S.:}{Rational curves and Extensions of Holomorphic
mappings.} Proc.\ Symp.\ Pure \ Math. {\bf52} Part 1, 93-104, (1991).

\entry{Iv-2}{Ivashkovich S.:}{Spherical shells as obstructions for the
extension of  holomorphic mappings.} The\ Journal\ of\ Geometric\ Analysis.
{\bf2}, N 4, 351-371, (1992).

\entry{Iv-3}{Ivashkovich S.:}{The Hartogs-type extension theorem for the
meromorphic maps into compact K\"ahler manifolds.} Invent.\ math. {\bf109}
 , 47-54, (1992).

\entry{Iv-4}{Ivashkovich S.:}{An example concerning extension and
separate analyticity properties of meromorphic mappings.} 
Amer. J. Math. {\bf 121}, 97-130 (1999).




 \entry{Ka-1}{Kato M.:}{Compact complex manifolds containing "global"
 spherical shells I.} Proc. Intl. Symp. Algebraic Geometry, Kyoto,
 45-84 (1977).
 
 \entry{Ka-2}{Kato M.:}{Examples of simply connected compact complex 3-folds.}
 Tokyo Journal of Mathematics {\bf5} 341-364 (1982).
 
 \entry{Ka-3}{Kato M.:}{Examples of simply connected compact complex 3-folds II.}
 Tokyo Journal of Mathematics {\bf9} 1-28 (1986).
 
\entry{Ka-4}{Kato M.:}{Examples on an Extension Problem of Holomorphic Maps
and Holomorphic 1-Dimensional Foliations.} Tokyo\ Journal\ Math. {\bf13}, n 1,
 139-146, (1990).


\entry{Kl}{Klimek M.:}{Pluripotential theory.} London.\ Math.\ Soc.\ Monographs
 ,\ New\ Series 6, (1991).



\entry{Lg}{Lelong P.:}{Plurisubharmonic Functions and Positive Differential
Forms. Gordon and Breach.} New-York. (1969) 78 p.

\entry{Lv-Sl}{Levenberg N., Slodkowski Z.}{Pseudoconcave pluripolar sets 
in $\cc^2$.} Math. Ann. {\bf312}, 429-443 (1998).
 
\entry{Lv}{Levi E.:}{Studii sui punti singolari essenziali delle funzioni
analitiche di due o pi\'u variabili complesse.} Annali di Mat. pura ed appl.
 {\bf17}, n 3, 61-87 (1910).

\entry{Mz}{Mazet P.:}{Un th\'eor\`eme d'image directe propre.} S\'eminaire 
P. Lelong, LNM {\bf 410}, 107-116 (1973).

\entry{Ni}{Nishino T.:}{Function Theory in Several Complex Variables.} AMS,  
Translation of Math. Monographs, v.193 (2000).

\entry{No}{Nori M.: }{The Schottky groups in higher dimension.} Contemporary 
Mathematics, {\bf58}, part I, 195-197 (1986).

\entry{Re}{Remmert R.:}{Holomorphe und meromorphe Abbildungen komplexer
R\"aume.} Math.\ Ann. {\bf133}, 328-370, (1957).



\entry{Sh-1}{Shiffman B.:}{Extension of Holomorphic Maps into Hermitian
Manifolds.} Math.\ Ann. {\bf194}, 249-258, (1971).

\entry{Sh-2}{Shiffman B.:}{Extension of positive line bundles and
meromorphic maps.} Invent. math. {\bf 15}, 332-347 (1976).

\entry{Sb}{Sibony N.:}{Quelques problemes de prolongement de courants en
analyse complexe.} Duke\ Math.\ J. {\bf52}, 157-197, (1985).

\entry{Si-1}{Siu Y.-T.:}{Every Stein subvariety admits a stein neighborhood.}
 Invent.\ Math. {\bf38}, N 1, 89-100, (1976).


\entry{Si-2}{Siu Y.-T.:}{Extension of meromorphic maps into K\"ahler
manifolds.} Ann.\ Math. {\bf102}, 421-462 (1975).

\entry{Si-3}{Siu Y.-T.:}{Analyticity of sets associated to Lelong numbers and the 
extension of closed positive currents.} Invent. Math. {\bf27}, 53-156 (1974).

\entry{Si-T}{Siu Y.-T., Trautmann G.:}{Gap-Sheaves and Extension of Coherent
Analytic Subsheaves.} Lect.\ Notes\ Math. {\bf172}, Springer-Verlag, (1971).




\entry{V}{Varouchas J.:}{Stabilit\'e de la classe des vari\'et\'es
K\"ahl\'eriennes
par certain morphismes propres.} Invent.\ math. {\bf77}, 117-127 (1984).

\bigskip\bigskip
Universit\'e de Lille-I

U.F.R. de Math\'ematiques

Villeneuve d'Ascq Cedex

59655 France

ivachkov@gat.univ-lille1.fr

\bigskip\bigskip
IAPMM Nat. Acad Sci. 

Naukova 3/b, 290053 Lviv

Ukraine

\end